\def\DateTime{6/Sep/2000, 11:50AM (JP)}
\def\Version{Version 2.1}

\def\yes{\if00}
\def\no{\if01}
\def\iftwelvept{\yes}
\def\ifusepdf{\no}
\def\ifpsfont{\yes}
\def\ifneedpagebreakforpic{\no}

\iftwelvept
\documentclass[leqno,12pt]{amsart}
\else
\documentclass[leqno]{amsart}
\fi
\usepackage{amssymb}
\usepackage{amscd}
\usepackage{latexsym}
\usepackage{verbatim}
\usepackage{pstcol}

\ifusepdf
\usepackage{hyperref}
\else\fi
\ifpsfont
\usepackage[T1]{fontenc}
\usepackage{times}
\else\fi

\input xy
\xyoption{all}

\iftwelvept
\else\fi

\theoremstyle{plain}
\newtheorem{Theorem}{Theorem}[section]
\newtheorem{Proposition}[Theorem]{Proposition}
\newtheorem{Lemma}[Theorem]{Lemma}
\newtheorem{Corollary}[Theorem]{Corollary}
\newtheorem{Claim}{Claim}[Theorem]

\theoremstyle{definition}

\newtheorem{Remark}[Theorem]{Remark}
\newtheorem{Example}[Theorem]{Example}

\newtheorem{Problem}[Theorem]{Problem}
\newtheorem{Question}[Theorem]{Question}

\renewcommand{\theTheorem}{\arabic{section}.\arabic{Theorem}}
\renewcommand{\theClaim}{\arabic{section}.\arabic{Theorem}.\arabic{Claim}}
\renewcommand{\theequation}{\arabic{section}.\arabic{Theorem}.\arabic{Claim}}

\def\rom{\textup}
\newcommand{\ZZ}{{\mathbb{Z}}}
\newcommand{\QQ}{{\mathbb{Q}}}

\newcommand{\PP}{{\mathbb{P}}}

\newcommand{\OO}{{\mathcal{O}}}
\newcommand{\MM}{{M}}
\newcommand{\bMM}{{\bar{M}}}
\newcommand{\codim}{\operatorname{codim}}
\newcommand{\dis}{\operatorname{dis}}

\newcommand{\Ker}{\operatorname{Ker}}
\newcommand{\Sing}{\operatorname{Sing}}
\newcommand{\ch}{\operatorname{char}}
\newcommand{\Coker}{\operatorname{Coker}}

\newcommand{\length}{\operatorname{length}}
\newcommand{\Sym}{\operatorname{Sym}}
\newcommand{\Supp}{\operatorname{Supp}}

\newcommand{\End}{\operatorname{\mathcal{E}\textsl{nd}}}
\newcommand{\Pic}{\operatorname{Pic}}
\newcommand{\Div}{\operatorname{Div}}
\newcommand{\Ref}{\operatorname{Ref}}
\newcommand{\Rat}{\operatorname{Rat}}

\newcommand{\PGL}{\operatorname{PGL}}
\newcommand{\Hilb}{\operatorname{Hilb}}
\newcommand{\zero}{\operatorname{div}}
\newcommand{\Proof}{{\sl Proof.}\quad}
\newcommand{\rank}{\operatorname{rk}}
\newcommand{\mult}{\operatorname{mult}}
\newcommand{\Cone}{\operatorname{Cone}}
\newcommand{\Nef}{\operatorname{Nef}}
\newcommand{\FNef}{\operatorname{FNef}}
\newcommand{\QED}{{\unskip\nobreak\hfil\penalty50\quad\null\nobreak\hfil
{$\Box$}\parfillskip0pt\finalhyphendemerits0\par\medskip}}
\newcommand{\rest}[2]{\left.{#1}\right\vert_{{#2}}}

\begin{document}

\title[Nef divisors on the moduli space of stable curves]%
{Nef divisors in codimension one \\
on the moduli space of stable curves}
\author{Atsushi Moriwaki}
\address{Department of Mathematics, Faculty of Science,
Kyoto University, Kyoto, 606-8502, Japan}
\email{moriwaki@kusm.kyoto-u.ac.jp}
\date{\DateTime, (\Version)}
\begin{abstract}
Let $\MM_g$ be the moduli space of 
smooth curves of genus $g \geq 3$, and
$\bMM_g$ the Deligne-Mumford compactification in terms
of stable curves.
Let $\bMM_g^{[1]}$ be an open set of $\bMM_g$ consisting of
stable curves of genus $g$
with one node at most.
In this paper, we determine the necessary and sufficient 
condition to guarantee that a $\QQ$-divisor $D$ on $\bMM_g$
is nef over $\bMM_g^{[1]}$, that is,
$(D \cdot C) \geq 0$ for all irreducible curves
$C$ on $\bMM_g$ with $C \cap \bMM_g^{[1]} \not= \emptyset$.
\end{abstract}


\maketitle

\tableofcontents

\section*{Introduction}
\renewcommand{\theTheorem}{\Alph{Theorem}}
Throughout this paper, we fix an algebraically closed field $k$, 
and every algebraic scheme is defined over $k$.
For simplicity, we assume that the characteristic of $k$ is zero
in this introduction.

Let $X$ be a normal complete variety and
$\mathcal{P}$ a certain kind of positivity of $\QQ$-line bundles on $X$ 
(e.g. ampleness, effectivity, bigness, etc).
A problem to describe the cone $\Cone(X; \mathcal{P})$ 
consisting of $\QQ$-line bundles with the positivity $\mathcal{P}$
is usually very hard and interesting.
In this paper, as positivity, we consider numerical effectivity
over a fixed open set.
Namely, let $U$ be a Zariski open set of $X$.
We say a $\QQ$-line bundle $L$ is {\em nef over $U$} if, 
for all irreducible curves $C$ with $C \cap U \not= \emptyset$,
$(L \cdot C) \geq 0$.
We define the relative nef cone $\Nef(X; U)$ over $U$
to be the cone of $\QQ$-line bundles on 
$X$ which are nef over $U$.

Let $g$ and $n$ be non-negative integers with $2g-2+n > 0$.
Let $\bMM_{g,n}$ (resp. $\MM_{g,n}$) denote the moduli space
of $n$-pointed stable curves (resp. $n$-pointed smooth curves) 
of genus $g$.
For a non-negative integer $t$,
an irreducible component of 
the closed subscheme consisting of curves with
at least $t$ nodes is called {\em a $t$-codimensional stratum}
of $\bMM_{g,n}$.
(For example, a $1$-codimensional stratum is a boundary component.)
We denote by $S^t(\bMM_{g,n})$ the set of all $t$-codimensional
strata of $\bMM_{g,n}$.
Let $\bMM_{g,n}^{[t]}$ be the open set of $\bMM_{g,n}$
obtained by subtracting all $(t+1)$-codimensional strata, i.e.,
$\bMM_{g,n}^{[t]}$ is the open set consisting of
curves with at most $t$ nodes.
(Note that $\bMM_{g,n}^{[0]} = \MM_{g,n}$.)
Here we consider the following problem:

\begin{Problem}
Describe the tower of relative nef cones
\[
\Nef(\bMM_{g,n};\MM_{g,n}) \supseteq 
\Nef(\bMM_{g,n};\bMM_{g,n}^{[1]}) \supseteq
\cdots \supseteq \Nef(\bMM_{g,n};\bMM_{g,n}^{[3g-3+n-1]}) = \Nef(\bMM_{g,n}).
\]
\end{Problem}

\noindent
We say a $\QQ$-divisor on $\bMM_{g,n}$ is {\em F-nef} 
if the intersection number with
every $1$-dimensional stratum is non-negative.
Let $\FNef(\bMM_{g,n})$ denote the cone consisting of F-nef $\QQ$-divisors.
Concerning the top $\Nef(\bMM_{g,n})$ of the tower,
it is conjectured in \cite{FICM}, \cite{KMCont} and \cite{GKM} that
$\FNef(\bMM_{g,n}) = \Nef(\bMM_{g,n})$.
In other words, the Mori cone of $\bMM_{g,n}$
is generated by $1$-dimensional strata,
which gives rise to a concrete description of $\Nef(\bMM_{g,n})$
(cf. \cite{FICM}, \cite{KMCont} and \cite{GKM}).
Moreover, it is closely related to the relative nef cone
$\Nef(\bMM_{g,n};\MM_{g,n})$. Actually, it was shown in \cite{GKM}
that if the weaker assertion 
$\FNef(\bMM_{g,n}) \subseteq \Nef(\bMM_{g,n};\MM_{g,n})$
holds for all $g,n$, then $\FNef(\bMM_{g,n}) = \Nef(\bMM_{g,n})$.
Further, as discussed in \cite{GKM},
$\bMM_{g,n}$ admits no interesting birational morphism to a projective variety.
However, we can expect the rich birational geometry on $\bMM_{g,n}$ in terms
of rational maps.
In this sense, to understand the tower of relative nef cones as above might
be a step toward this natural problem.

\medskip
We assume that $g \geq 3$ and $n=0$.
Let $\lambda$ be the Hodge class on $\bMM_g$, and 
$\delta_{irr}, \delta_1, \ldots, \delta_{[g/2]}$
the classes of irreducible components of the boundary 
$\bMM_g \setminus \MM_g$.
Let $\mu$ be a divisor on $\bMM_g$ given by
\[
\mu = (8g+4) \lambda - g \delta_{irr} - \sum_{i=1}^{[g/2]}4i(g-i)
\delta_i.
\]
In the paper \cite{MoRB}, we 
proved that $\Nef(\bMM_g; \MM_g)$ is the convex hull
spanned by $\mu, \delta_{irr}, \delta_1, \ldots, \delta_{[g/2]}$,
that is, 
\[
\Nef(\bMM_g; \MM_g) =
\QQ_{+} \mu + \QQ_{+} \delta_{irr} + \sum_{i=1}^{[g/2]} \QQ_{+} \delta_i,
\]
where $\QQ_{+} = \{ x \in \QQ \mid x \geq 0 \}$.
The main result of this paper is to determine 
$\Nef(\bMM_g;\bMM_g^{[1]})$.

\begin{Theorem}[cf. Theorem~\ref{thm:cone:mg:codim:two}]
\label{thm:cone:mg:codim:two:intro}
A $\QQ$-divisor $a \mu + b_{irr} \delta_{irr} +
\sum_{i=1}^{[g/2]} b_{i}
\delta_{i}$ on $\bMM_{g}$
is nef over $\bMM_{g}^{[1]}$ if and only if
the following system of inequalities hold:
\begin{align*}
& a \geq \max \left\{ \frac{b_{i}}{4i(g-i)} \mid 
i = 1, \ldots, [g/2] \right\}, \\
& B_0 \geq B_1 \geq B_2 \geq \cdots \geq B_{[g/2]}, \\
& B^*_{[g/2]} \geq \cdots \geq B^*_2 \geq B^*_1 \geq B^*_0,
\end{align*}
where $B_0$, $B_0^*$, $B_i$ and $B_i^*$ 
\rom{(}$i=1, \ldots, [g/2]$\rom{)} are given by
\[
B_0 = 4b_{irr},\quad
B^*_0 = \frac{4b_{irr}}{g(2g-1)},\quad
B_i = \frac{b_i}{i(2i+1)}\quad\text{and}\quad
B^*_i = \frac{b_i}{(g-i)(2(g-i)+1)}.
\]
\end{Theorem}

An interesting point is that the above theorem shows us 
that $\mu$ is not only nef over $\MM_{g}$ but also 
nef over $\bMM_g^{[1]}$.
Moreover, the theorem tells us that every nef $\QQ$-divisor
over $\bMM_g^{[1]}$ can be obtained in the following way.
Namely, we first fix a non-negative rational number
$b_{irr}$, and take $b_1$ with
\[
\frac{4(g-1)b_{irr}}{g} \leq b_1 \leq 12 b_{irr}.
\]
Further, we choose $b_2, \ldots, b_{[g/2]}$ inductively by using
\[
\frac{(g-1-i)(2(g-i)-1)}{(g-i)(2(g-i)+1)}b_i \leq b_{i+1} \leq
\frac{(i+1)(2i+3)}{i(2i+1)}b_i.
\]
Finally, we take $a$ with
\[
a \geq \max \left\{ \frac{b_{i}}{4i(g-i)} \mid 
i = 1, \ldots, [g/2] \right\}.
\]
Then, a $\QQ$-divisor given by
$a \mu + b_{irr}\delta_{irr} + \sum_{i=1}^{[g/2]}b_i \delta_i$
is nef over $\bMM_g^{[1]}$.
Further, as corollaries of the above theorem, we have the following.

\begin{Corollary}[cf. Corollary~\ref{cor:cone:mg:codim:two}]
\label{cor:cone:mg:codim:two:intro:1}
For an irreducible component $\Delta$ of the boundary $\bMM_{g} \setminus
\MM_g$, let $\widetilde{\Delta}$ be the normalization of
$\Delta$, and $\rho_{\Delta} : \widetilde{\Delta} \to \bMM_{g}$
the induced morphism.
Then, a $\QQ$-divisor $D$ on $\bMM_g$ is nef over $\bMM_g^{[1]}$ if and only if
the following are satisfied:
\begin{enumerate}
\renewcommand{\labelenumi}{(\arabic{enumi})}
\item
$D$ is weakly positive at any points of $\MM_g$.

\item
For every boundary component $\Delta$,
$\rho_{\Delta}^*(D)$ is weakly positive at any points of 
$\rho_{\Delta}^{-1}(\bMM_g^{[1]})$
\end{enumerate}
For the definition of weak positivity, see \rom{\S\ref{subsec:pos:weil:div}}.
\end{Corollary}

\begin{Corollary}[cf. Corollary~\ref{cor:cone:mg:codim:two:another}]
\label{cor:cone:mg:codim:two:intro:2}
With notation as above,
if $\rho_{\Delta}^*(D)$ is nef over $\rho_{\Delta}^{-1}(\bMM_g^{[1]})$
for every boundary component $\Delta$, then $D$ is nef over $\bMM_{g}^{[1]}$.
In particular, the Mori cone of $\bMM_g$ is the convex hull spanned
by curves lying on the boundary $\bMM_g \setminus \MM_g$,
which gives rise to a special case of \cite[Proposition~3.1]{GKM}.
\end{Corollary}

Let us go back to the general situation.
Similarly, for $\Delta \in S^l(\bMM_{g,n})$, 
let $\widetilde{\Delta}$ be the normalization of
$\Delta$, and $\rho_{\Delta} : \widetilde{\Delta} \to \bMM_{g,n}$
the induced morphism.
Inspired by the above corollaries, we have the following questions:

\begin{Question}
\label{question:nef:wp:intro}
For a non-negative integer $t$, if a $\QQ$-divisor $D$ on $\bMM_{g,n}$
is nef over $\bMM_{g,n}^{[t]}$, then 
is $\rho_{\Delta}^*(D)$ weakly positive at any points of
$\rho_{\Delta}^{-1}(\bMM_{g,n}^{[l]})$
for all $0 \leq l \leq t$ and all $\Delta \in S^l(\bMM_{g,n})$?
More strongly, if
$D$ is nef over $\bMM_{g,n}^{[t]}$, then
is $D$ weakly positive at any points of $\bMM_{g,n}^{[t]}$?
\end{Question}

\begin{Question}
\label{question:nef:in:small:nef:in:bigintro}
Fix an integer $t$ with $0 \leq t \leq 3g-3+n-1$.
If $\rho_{\Delta}^*(D)$ is nef over $\rho_{\Delta}^{-1}(\bMM_{g,n}^{[t]})$
for all $\Delta \in S^t(\bMM_{g,n})$, then is $D$ is nef over $\bMM_{g,n}^{[t]}$?
\end{Question}

\noindent
In the case $t = 3g-3+n-1$, the above question is nothing more than
asking $\FNef(\bMM_{g,n}) = \Nef(\bMM_{g,n})$.

\medskip
In order to get the above theorem,
we need a certain kind of slope inequalities on 
the moduli space of $n$-pointed stable curves.
The $\QQ$-line bundles $\lambda$ and $\psi_1, \ldots, \psi_n$ on
$\bMM_{g,n}$ are defined as follows:
Let $\pi : \bMM_{g,n+1} \to \bMM_{g,n}$ be the universal curve
of $\bMM_{g,n}$, and $s_1, \ldots, s_n : \bMM_{g,n} \to \bMM_{g,n+1}$
the sections of $\pi$ arising from the $n$-points of $\bMM_{g,n}$.
Then, $\lambda = \det(\pi_* (\omega_{\bMM_{g,n+1}/\bMM_{g,n}}))$ and
$\psi_i = s_i^*(\omega_{\bMM_{g,n+1}/\bMM_{g,n}})$ for $i=1, \ldots, n$.
Here we set
\begin{align*}
[n] & = \{ 1, \ldots, n\}\quad (\text{note that $[0] = \emptyset$}), \\
\Upsilon_{g,n} & = \{ (i, I) \mid \text{$i \in \ZZ$, $0 \leq i \leq g$
and $I \subseteq [n]$} \} \setminus 
\{ (0, \emptyset), (0, \{1\}), \ldots, (0, \{n\})  \}, \\
\overline{\Upsilon}_{g,n} & = \{ \{(i,I), (j,J)\} \mid 
(i,I), (j,J) \in \Upsilon_{g,n}, i+j=g, I \cap J = \emptyset, I\cup J = [n] \}.
\end{align*}
Moreover, for a finite set $S$, we denote the number of it by $\vert S \vert$.
The boundary
$\bMM_{g,n} \setminus \MM_{g,n}$ has 
the following irreducible decomposition:
\[
 \bMM_{g,n} \setminus \MM_{g,n} = \Delta_{irr} \cup 
 \bigcup_{\{ (i,I), (j,J) \} \in \overline{\Upsilon}_{g,n}} 
 \Delta_{\{ (i, I), (j, J) \}}.
\]
A general point of $\Delta_{irr}$ represents an $n$-pointed
irreducible stable curve with one node.
A general point of $\Delta_{\{(i,I), (j,J)\}}$ represents an $n$-pointed
stable curve consisting of an $\vert I \vert$-pointed smooth curve 
$C_1$ of genus $i$ and
a $\vert J \vert$-pointed smooth curve $C_2$ of genus $j$ 
meeting transversally
at one point, where $\vert I \vert$-points on $C_1$ 
(resp. $\vert J \vert$-points on $C_2$)
arise from $\{ s_t \}_{t \in I}$ (resp. $\{ s_l \}_{l \in J}$).
Let $\delta_{irr}$ and $\delta_{\{(i,I),(j,J)\}}$ be the classes of
$\Delta_{irr}$ and $\Delta_{\{(i,I),(j,J)\}}$ in 
$\Pic(\bMM_{g,n}) \otimes \QQ$ respectively.
For a subset $L$ of $[n]$, we define a $\QQ$-divisor 
$\theta_L$ on $\bMM_{g,n}$ to be
\[
\theta_L = 
4(g-1+\vert L \vert)(g-1) \sum_{t \in L} \psi_t - 12 \vert L \vert^2 \lambda 
+ \vert L \vert^2  \delta_{irr} 
- \sum_{\upsilon \in \overline{\Upsilon}_{g,n}} 4\gamma_L(\upsilon)
\delta_{\upsilon},
\]
where $\gamma_L : \overline{\Upsilon}_{g,n} \to \ZZ$ is given by
\[
\gamma_L\left(\{ (i,I), (j,J) \}\right) =
\left( \det \begin{pmatrix} i & \vert L \cap I \vert \\ j & 
\vert L \cap J \vert \end{pmatrix}
+ \vert L \cap I \vert \right)
\left( \det \begin{pmatrix} i & \vert L \cap I \vert \\ j & 
\vert L \cap J \vert \end{pmatrix} 
- \vert L \cap J \vert \right).
\]
Then, we have the following.

\begin{Theorem}[cf. Theorem~\ref{thm:slope:inq:m:g:n}]
\label{thm:slope:inq:m:g:n:intro}
For any subset $L$ of $[n]$, the divisor $\theta_L$
is weakly positive at any points of $\MM_{g,n}$.
In particular, it is nef over $\MM_{g,n}$.
\end{Theorem}

\noindent
We remark that
R. Hain has already announced the above inequality in the case where $n=1$.
(For details, see \cite{Hain}.)
Theorem~\ref{thm:slope:inq:m:g:n:intro} is a generalization of his inequality.

\medskip
Here we assume that $g \geq 2$.
First note that
\[
\mu = (8g+4)\lambda - g \delta_{irr} - 
\sum_{\{ (i,I), (j,J) \} \in \overline{\Upsilon}_{g,n}} 
4ij \delta_{\{ (i, I), (j, J) \}}
\]
is nef over $\MM_{g,n}$.
Thus, as a consequence of Theorem~\ref{thm:slope:inq:m:g:n:intro},
we can see that
\[
\QQ_{+} \mu + \sum_{L \subseteq [n]} \QQ_{+} \theta_L + \QQ_{+} \delta_{irr}
+ \sum_{ \upsilon \in \overline{\Upsilon}_{g,n}} 
\QQ_{+} \delta_{\upsilon}
\subseteq \Nef(\bMM_{g,n}; \MM_{g,n}),
\]
so that we may ask the following question:

\begin{Question}
Is $\Nef(\bMM_{g,n}; \MM_{g,n})$ the convex hull spanned by $\QQ$-divisors
$\mu$, $\theta_L$ ($\forall L \subseteq [n]$), $\delta_{irr}$ and
$\delta_{\upsilon}$ ($\forall \upsilon \in 
\overline{\Upsilon}_{g,n}$).
\end{Question}

\noindent
Corollary~\ref{cor:slope:cone:mg1} and 
Corollary~\ref{cor:slope:cone:mg2} are partial answers for the above question.
If the above question is true, then it gives an affirmative answer of
Question~\ref{question:nef:wp:intro} for $t=0$.

\medskip
Finally, we would like to give hearty thanks to Prof. Hain and
Prof. Keel for their useful comments for this paper.

\renewcommand{\theTheorem}{\arabic{section}.\arabic{subsection}.\arabic{Theorem}}
\renewcommand{\theClaim}{\arabic{section}.\arabic{subsection}.\arabic{Theorem}.\arabic{Claim}}
\renewcommand{\theequation}{\arabic{section}.\arabic{subsection}.\arabic{Theorem}.\arabic{Claim}}

\section{Notations, conventions, terminology and preliminaries}

Throughout this paper, we fix an algebraically closed field $k$,
and every algebraic scheme is defined over $k$.

\subsection{The positivity of Weil divisors}
\label{subsec:pos:weil:div}
\setcounter{Theorem}{0}
Let $X$ be a normal variety.
Let denote $Z^1(X)$ (resp. $\Div(X)$) the group of Weil divisors 
(resp. Cartier divisors) on $X$, and
$\sim$ the linear equivalence on $Z^1(X)$.
We set $A^1(X) = Z^1(X)/\!\!\!\sim$ and $\Pic(X) = \Div(X)/\!\!\!\sim$.
Note that $\Pic(X)$ is canonically isomorphic to the
Picard group (the group of isomorphism classes of line bundles).
Moreover, we denote by $\Ref(X)$ the set of isomorphism classes
of reflexive sheaves of rank $1$ on $X$.
For a Weil divisor $D$, the sheaf $\OO_X(D)$ is given by
\[
\OO_X(D)(U) = \{ \phi \in \Rat(X)^{\times} \mid
\text{$(\phi) + D$ is effective over $U$} \} \cup \{ 0 \}
\]
for each Zariski open set $U$ of $X$.
Then, we can see $\OO_X(D) \in \Ref(X)$.
Conversely, let $L$ be a reflexive sheaf of rank $1$ on $X$.
For a non-zero rational section $s$ of $L$,
$\zero(s)$ is defined as follows:
Let $X_0$ be the maximal Zariski open set of $X$ over which 
$L$ is locally free.
Note that $\codim(X \setminus X_0) \geq 2$.
Then, $\zero(s) \in Z^1(X)$ is defined by the Zariski closure of
$\zero(\rest{s}{X_0})$.
By our definition, we can see that $\OO_X(\zero(s)) \simeq L$.
Thus, 
the correspondence $D \mapsto \OO_X(D)$ gives rise to an
isomorphism $A^1(X) \simeq \Ref(X)$.
Here we remark that
if $x \not\in \Supp(\zero(s))$, then $L$ is free at $x$ because
$\OO_X(\zero(s))_x = \OO_{X,x}$ for $x \not\in \Supp(\zero(s))$.

An element of $Z^1(X) \otimes \QQ$ (resp. $\Div(X) \otimes \QQ$)
is celled a {\em $\QQ$-divisor} (resp. {\em $\QQ$-Cartier divisor}).
For $\QQ$-divisors $D_1$ and $D_2$, we say $D_1$ is 
$\QQ$-linearly equivalent to
$D_2$, denoted by $D_1 \sim_{\QQ} D_2$, if
there is a positive integer $n$ such that
$nD_1, nD_2 \in Z^1(X)$ and $nD_1 \sim nD_2$, i.e.,
$D_1$ coincides with $D_2$ in $A^1(X) \otimes \QQ$.

Fix a subset $S$ of $X$.
For $D \in Z^1(X) \otimes \QQ$, we say $D$ is {\em semi-ample over $S$}
if, for any $s \in S$, there is an effective $\QQ$-divisor $E$ on $X$ with
$s \not\in \Supp(E)$ and $D \sim_{\QQ} E$.
Moreover, $D$ is said to be 
{\em weakly positive over $S$}
if there are $\QQ$-divisors $Z_1, \ldots, Z_{l}$, a sequence
$\{ D_m \}_{m=1}^{\infty}$ of $\QQ$-divisors, and
sequences 
$\{ a_{1, m} \}_{m=1}^{\infty}, \ldots, \{a_{l, m} \}_{m=1}^{\infty}$
of rational numbers such that 
\begin{enumerate}
\renewcommand{\labelenumi}{(\arabic{enumi})}
\item
$l$ does not depend on $m$,

\item
$D_m$ is semi-ample over $S$ for all $m \gg 0$,

\item
$D \sim_{\QQ} D_m + \sum_{i=1}^{l} a_{i, m} Z_i$ for all $m \gg 0$, and

\item
${\displaystyle \lim_{m\to\infty} a_{i, m} = 0}$ for all $i=1, \ldots, l$.
\end{enumerate}
In the above definition, if $D$, $D_m$ and $Z_i$'s are
$\QQ$-Cartier divisors, then $D$ is said to be
{\em weakly positive over $S$ in terms of Cartier divisors}
(for short, C-weakly positive over $S$).
Further, if $D$ is semi-ample over $\{ x \}$ for some $x \in X$,
then we say $D$ is semi-ample at $x$.
Similarly, we define
the weak positivity of $D$ at $x$ and the C-weak positivity of $D$ at $x$.
We remark that weak positivity in \cite{MoRB} is 
nothing more than C-weak positivity.
Moreover, note that if a $\QQ$-divisor $D$ is semi-ample at $x$,
then $D$ is a $\QQ$-Cartier divisor around $x$, i.e.,
there is a Zariski open set $U$ of $X$ such that
$x \in U$ and $\rest{D}{U}$ is a $\QQ$-Cartier divisor on $U$.

A normal variety $X$ is said to be {\em $\QQ$-factorial} if 
$Z^1(X) \otimes \QQ = \Div(X) \otimes \QQ$, i.e.,
any Weil divisors are $\QQ$-Cartier divisors.
It is well known that if $Y \to X$ is a finite and surjective morphism
of normal varieties and $Y$ is $\QQ$-factorial, then
$X$ is also $\QQ$-factorial (cf. \cite[Lemma~5.16]{KMBirat}).
Thus the moduli space $\bMM_{g,n}$ of $n$-pointed stable curves
of genus $g$
is $\QQ$-factorial because $\bMM_{g,n}$ is
an orbifold.
If $X$ is $\QQ$-factorial, then
the weak positivity of $D$ over $S$ coincides with
the C-weak positivity of $D$ over $S$.

We assume that $X$ is complete and $D$ is a $\QQ$-Cartier divisor.
We say $D$ is {\em nef over $S$} if $(D \cdot C) \geq 0$ 
for any complete irreducible curves $C$ with $S \cap C \not= \emptyset$.
Moreover, for a point $x$ of $X$,
we say {\em $D$ is nef at $x$} if $D$ is nef over $\{ x \}$.
Note that
\[
\text{``$D$ is semi-ample at $x$''} \quad \Longrightarrow \quad
\text{``$D$ is C-weakly positive at $x$''} 
\quad \Longrightarrow \quad
\text{``$D$ is nef at $x$''}
\]

\begin{Lemma}[$\ch(k) \geq 0$]
\label{lem:semi:finite:points}
Let $D$ be a $\QQ$-divisor on $X$, and $x_1, \ldots, x_n \in X$.
If $D$ is semi-ample at $x_i$ for each $i$, then
there is an effective $\QQ$-divisor $E$ on $X$ such that
$E \sim_{\QQ} D$ and $x_i \not\in \Supp(E)$ for all $i$.
\end{Lemma}

\Proof
By our assumption,
there is an effective $\QQ$-divisor $E_i$ on $X$ such that
$E_i \sim_{\QQ} D$ and $x_i \not\in \Supp(E_i)$.
Take a sufficiently large integer $m$ such that
$mD, mE_1, \ldots, mE_n \in Z^1(X)$ and
$mD \sim mE_i$ for all $i$.
Thus, there is a section $s_i$ of
$H^0(X, \OO_X(mD))$ with $\zero(s_i) = mE_i$.
Here since $x_i \not\in \Supp(mE_i)$ and $\OO_X(mD) \simeq \OO_X(mE_i)$,
we can see that $\OO_X(mD)$ is free at each $x_i$.

For $\alpha = (\alpha_1, \ldots\, \alpha_n) \in k^n$,
we set $s_{\alpha} = \alpha_1 s_1 + \cdots + \alpha_n s_n \in 
H^0(X, \OO_X(mD))$.
Further, we set $V_i = \{ \alpha \in k^n \mid s_{\alpha}(x_i) = 0 \}$.
Then, $\dim V_i= n-1$ for all $i$. Thus, since
$\#(k) = \infty$, there is $\alpha \in k^n$
with $\alpha \not\in V_1 \cup \cdots \cup V_r$, i.e.,
$s_{\alpha}(x_i) \not= 0$ for all $i$.
Let us consider a divisor $E = \operatorname{div}(s_{\alpha})$. 
Then, $E \sim m D$
and $x_i \not\in \Supp(E)$ for all $i$.
\QED

\begin{Proposition}[$\ch(k) \geq 0$]
\label{prop:samp:pamp:push}
Let $\pi : X \to Y$ be a surjective, proper and generically finite morphism
of normal varieties.
Let $D$ be a $\QQ$-divisor on $X$ and $S$ a subset of $Y$ such that
$\pi^{-1}(S)$ is finite.
Then, we have the following.
\begin{enumerate}
\renewcommand{\labelenumi}{(\arabic{enumi})}
\item
If $D$ is semi-ample over $\pi^{-1}(S)$, then
$\pi_*(D)$ is semi-ample over $S$.

\item
If $D$ is weakly positive over $\pi^{-1}(S)$, then
$\pi_*(D)$ is weakly positive over $S$.
\end{enumerate}
\end{Proposition}

\Proof
(1) By Lemma~\ref{lem:semi:finite:points}, there is an effective divisor
$E$ on $X$ such that $E \sim_{\QQ} D$ and $s' \not\in \Supp(E)$ for all 
$s' \in \pi^{-1}(S)$. Then, $\pi_*(E) \sim_{\QQ} \pi_*(D)$ and
$s \not\in \pi(\Supp(E)) = \Supp(\pi_*(E))$ for all $s \in S$.

(2) This is a consequence of (1).
\QED

\begin{Proposition}[$\ch(k) \geq 0$]
\label{prop:samp:pamp:pullback}
Let $\pi : X \to Y$ be a surjective, proper morphism
of normal varieties. We assume that $Y$ is $\QQ$-factorial.
Let $D$ be a $\QQ$-divisor on $Y$, and $S$ a subset of $Y$.
Then, we have the following.
\begin{enumerate}
\renewcommand{\labelenumi}{(\arabic{enumi})}
\item
If $D$ is semi-ample over $S$, then
$\pi^*(D)$ is semi-ample over $f^{-1}(S)$.

\item
If $D$ is weakly positive over $S$, then
$\pi^*(D)$ is C-weakly positive over $S$.
\end{enumerate}
\end{Proposition}

\Proof
(1) Let $s'$ be a point in $\pi^{-1}(S)$.
Then, there is an effective $\QQ$-divisor $E$ on $Y$ with
$D \sim_{\QQ} E$ and $\pi(s') \not\in \Supp(E)$.
Thus, $\pi^*(D) \sim_{\QQ} \pi^*(E)$ and $s' \not\in \Supp(\pi^*(E))$.
Therefore, $\pi^*(D)$ is semiample over $\pi^{-1}(S)$.

(2) This is a consequence of (1).
\QED

\begin{Lemma}[$\ch(k) \geq 0$]
\label{lem:nef:over:product}
Let $X$ and $Y$ be complete varieties, and 
let $D$ and $E$ be $\QQ$-Cartier divisors on $X$ and $Y$ respectively.
Let $p : X \times Y$ and $q : X \times Y \to Y$ be the projections
to the first factor and the second factor respectively.
For $(x, y) \in X \times Y$,
$p^*(D) + q^*(E)$ is nef at $(x, y)$ if and only if
$D$ and $E$ are nef at $x$ and $y$ respectively.
\end{Lemma}

\Proof
First we assume that $p^*(D) + q^*(E)$ is nef at $(x, y)$.
Let $C$ be a complete irreducible curve on $X$ with $x \in C$.
Then, $C_y = C \times \{ y \}$ is a complete curve on
$X \times Y$ with $(x, y) \in C_y$.
Moreover, $(p^*(D) + q^*(E) \cdot C_y) = (D \cdot C)$.
Thus, $(D \cdot C) \geq 0$, which says us that
$D$ is nef at $x$. In the same way, we can see that
$E$ is nef at $y$.

Next we assume that $D$ and $E$ are nef at $x$ and $y$ respectively.
In order to see that $p^*(D) + q^*(E)$ is nef at $(x, y)$,
it is sufficient to check that $(p^*(D) \cdot C) \geq 0$ and
$(q^*(E) \cdot C) \geq 0$
for any complete irreducible curves $C$ on $X \times Y$ with
$(x, y) \in C$.
Here, $p(C)$ is either $\{ x \}$, or a complete irreducible curve
passing through $x$. Thus, by virtue of the projection formula,
$(p^*(D) \cdot C) \geq 0$. In the same way, $(q^*(E) \cdot C) \geq 0$.
\QED

\subsection{The first Chern class of coherent sheaves}
\setcounter{Theorem}{0}
Let $X$ be a normal variety, and
$F$ a coherent $\OO_X$-module on $X$.
Here we define $c_1(F) \in A^1(X)$ in the following way.

Case 1. $F$ is a torsion sheaf. In this case,
we set
\[
D = \sum_{\substack{P \in X, \\ \operatorname{depth}(P) = 1}} 
\operatorname{lenght}(F_P) \overline{ \{ P \} },
\]
where $\overline{ \{ P \} }$ is the Zariski closure of
$\{ P \}$ in $X$.
Then, $c_1(F)$ is defined by the class of $D$.

Case 2. $F$ is a torsion free sheaf.
Let $r$ be the rank of $F$. Then, $(\bigwedge^{r} F)^{\vee\vee}$
is a reflexive sheaf of rank $1$, where
${}^{\vee\vee}$ means the double dual of sheaves.
Thus, we define $c_1(F)$ to be the class of $(\bigwedge^{r} F)^{\vee\vee}$.

Case 3. $F$ is general. Let $T$ be the torsion part of $F$.
Then, $c_1(F) = c_1(T) + c_1(F/T)$.

\medskip
Note that if $0 \to F_1 \to F_2 \to F_3 \to 0$ is an exact sequence of
coherent $\OO_X$-modules, then $c_1(F_2) = c_1(F_1) + c_1(F_3)$.
Moreover, let $L$ be a reflexive sheaf of rank $1$ on $X$, and
$s$ a non-zero section of $L$. Then
\[
c_1(L) =
c_1\left(  \Coker(\OO_X \overset{\times s}{\longrightarrow} L) \right) 
= \text{the class of $\zero(s)$}.
\]

\begin{Proposition}[$\ch(k) \geq 0$]
\label{prop:samp:det}
Let $X$ be a normal algebraic variety,
$F$ a coherent $\OO_X$-module, and $x$ a point of $X$.
If $F$ is generated by global sections at $x$ and
$F$ is free at $x$, then $c_1(F)$ is
semi-ample at $x$.
\end{Proposition}

\Proof
Let $T$ be the torsion part of $F$. Then,
$c_1(F) = c_1(F/T) + c_1(T)$.
Here since $F$ is free at $x$, 
$c_1(T)$ is semi-ample at $x$.
Moreover, it is easy to see that
$F/T$ is generated by global sections at $x$.
Therefore, to prove our proposition, 
we may assume that $F$ is a torsion free sheaf.

Let $r$ be the rank of $F$ and
$\kappa(x)$ the residue field of $x$.
Then, by our assumption,
there are sections $s_1, \ldots, s_r$ of $F$ such that
$\{ s_i(x) \}$ forms a basis of $F \otimes \kappa(x)$.
Since we can view $s_i$ as an injection $s_i : \OO_X \to F$,
$s = s_1 \wedge \cdots \wedge s_r$ gives rise to
an injection $s : \OO_X \to \left( \bigwedge^r F \right)^{\vee\vee}$, which is
bijective at $x$.
Thus, $x \not\in \zero(s)$.
\QED

\subsection{The discriminant divisor of vector bundles}
\setcounter{Theorem}{0}
Let $f : X \to Y$ be a proper surjective morphism of algebraic varieties
of the relative dimension one, and let $E$ be a locally
free sheaf on $X$. We define the {\em discriminant divisor of $E$
with respect to $f$} to be
\[
 \dis_{X/Y}(E) = f_* \left( 2 \rank(E)c_2(E) - (\rank(E) -1 )c_1(E)^2
\right).
\] 

\begin{Lemma}[$\ch(k) \geq 0$]
\label{lem:comp:dis}
Let $f : X \to Y$ be a flat, surjective and
projective morphism of varieties with $\dim f = 1$.
Let $E$ be a vector bundle of rank $r$ on $X$.
Then, we have the following.
\begin{enumerate}
\renewcommand{\labelenumi}{(\arabic{enumi})}
\item
$\dis_{X/Y}(E)$ is a Cartier divisor.

\item
Let $u : Y' \to Y$ be a morphism of varieties, and let
\[
\begin{CD}
X @<{u'}<< X \times_Y Y' \\
@V{f}VV @VV{f'}V \\
Y @<{u}<< Y'
\end{CD}
\]
be the induced diagram of the fiber product.
If $X \times_Y Y'$ is integral, then $\dis_{X \times_Y Y'/Y'}({u'}^*(E)) 
= u^*(\dis_{X/Y}(E))$.
\end{enumerate}
\end{Lemma}

\Proof
(1)
We set $F = \End(E)$.
Let $p : P = \PP(F) \to X$ be the projective bundle of $F$, 
and $\OO_{P}(1)$ the
tautological line bundle on $P$.
Let $g : P \to Y$ be the composition of $P \overset{p}{\longrightarrow} X
\overset{f}{\longrightarrow} Y$.
Then, since
\[
p_*(c_1(\OO_P(1))^{r^2+1}) = -c_2(F) = -(2r c_2(E) - (r-1)c_1(E)^2),
\]
we have $g_*(c_1(\OO_P(1))^{r^2+1}) = -\dis_{X/Y}(E)$.
Thus, $\dis_{X/Y}(E) = 
-c_1\left(\langle \OO_P(1)^{\cdot r^2 +1} \rangle(P/Y)\right)$,
where
\[
\langle \ , \cdots, \ \rangle(P/Y) : 
\overbrace{\Pic(P) \times \cdots \times \Pic(P)}^{\dim g + 1}
\to \Pic(Y)
\]
is Deligne's pairing for
the flat morphism $g : P \to Y$.
Therefore, $\dis_{X/Y}(E)$ is a Cartier divisor.

(2)
This follows from the compatibility of 
Deligne's pairing by base changes.
\QED

\begin{Remark}
In (2) of Lemma~\ref{lem:comp:dis},
$X \times_Y Y'$ is integral if the generic fiber of
$X \times_Y Y' \to Y'$ is integral
by virtue of \cite[Lemma~4.2]{MoCD}.
\end{Remark}

\subsection{The moduli space of $T$-pointed stable curves of genus $g$}
\label{subsec:n:point:stable:curve}
\setcounter{Theorem}{0}
Let $g$ be a non-negative integer and
$T$ a finite set with $2g-2 + \vert T \vert > 0$,
where $\vert T \vert$ is the number of $T$.
Recall that $[n] = \{1, \ldots, n\}$ and $[0] = \emptyset$.
Usually, we use $[n]$ as $T$.
Let $\bMM_{g,T}$ (resp. $\MM_{g,T}$) denote the moduli space
of $T$-pointed stable curves (resp. $T$-pointed smooth curves)
of genus $g$,
namely, $\bMM_{g,T}$ (resp. $\MM_{g,T}$) is the moduli space
of $\vert T \vert$-pointed stable curves (resp. $\vert T \vert$-pointed
smooth curves) of genus $g$, whose marked points are labeled by
the index set $T$.

Roughly speaking,
the $\QQ$-line bundles $\lambda$ and $\{ \psi_t \}_{t \in T}$ on
$\bMM_{g,T}$ are defined as follows:
Let $\pi : \mathcal{C} \to \bMM_{g,T}$ be the universal curve
of $\bMM_{g,T}$, and $s_t : \bMM_{g,T} \to \mathcal{C}$ ($t \in T$)
the sections of $\pi$ arising from the $T$-points of $\bMM_{g,T}$.
Then, $\lambda = \det(\pi_* (\omega_{\mathcal{C}/\bMM_{g,T}}))$ and
$\psi_t = s_t^*(\omega_{\mathcal{C}/\bMM_{g,T}})$ for $t \in T$.

For $x \in \bMM_{g,T}$, let denote $C_x$ the nodal curve corresponding to $x$
(here we forget the $T$-points).
Let $S^l(\bMM_{g,T})$ be the set of all irreducible components of
the closed set
\[
\{ x \in \bMM_{g,T} \mid \#(\Sing(C_x)) \geq l \}.
\]
Then, every element of $S^l(\bMM_{g,T})$ is of codimension $l$,
so that it is called an {\em $l$-codimensional stratum}
of $\bMM_{g,T}$. Note that $\bMM_{g,T} \setminus \MM_{g,T}$ is a normal
crossing divisor in the sense of orbifolds. Thus
the normalization of an element of $S^l(\bMM_{g,T})$
is $\QQ$-factorial.
Moreover, we set 
\[
\bMM_{g,T}^{[l]} = 
\bMM_{g,T} \setminus \bigcup_{\Delta \in S^{l+1}(\bMM_{g,T})} \Delta,
\]
i.e., $\bMM_{g,T}^{[l]} = \{ x \in \bMM_{g,T} \mid \#(\Sing(C_x)) \leq l \}$.
Note that $\bMM_{g,T}^{[0]} = \MM_{g,T}$.

To describe the boundary of $\bMM_{g,T}$,
we set
\begin{align*}
\Upsilon_{g,T} & = \{ (i, I) \mid 
\text{$i \in \ZZ$, $0 \leq i \leq g$ and $I \subseteq T$} \} \setminus 
\left( \{ (0, \emptyset) \} \cup \{ (0, \{t\}) \}_{t \in T} \right), \\
\overline{\Upsilon}_{g,T} & = \{ \{(i,I), (j,J)\} \mid 
(i,I), (j,J) \in \Upsilon_{g,T}, i+j=g, I \cap J = \emptyset, I\cup J = T \}.
\end{align*}
Then, the boundary $\Delta = \bMM_{g,T} \setminus \MM_{g,T}$
has the following irreducible decomposition:
\[
 \Delta = \Delta_{irr} \cup 
 \bigcup_{\{ (i,I), (j,J) \} \in \overline{\Upsilon}_{g,T}} 
 \Delta_{\{(i,I), (j,J)\}}.
\]
A general point of $\Delta_{irr}$ represents a $T$-pointed
irreducible stable curve with one node.
A general point of $\Delta_{\{(i,I), (j,J)\}}$ represents a $T$-pointed
stable curve consisting of an $I$-pointed smooth curve of genus $i$ and
a $J$-pointed smooth curve of genus $j$ meeting transversally
at one point.
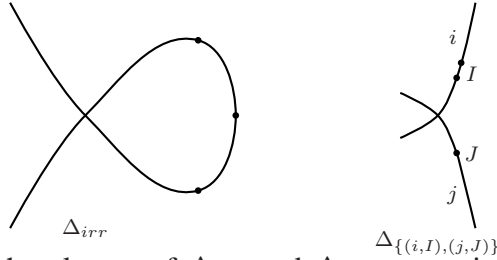
\[
\begin{pspicture}(0,0)(3,3)
\pscurve(0,0)(1,1.5)(2.5,2.5)(3,1.5)(2.5,0.5)(1,1.5)(0,3)
\psdots(2.5,2.5)(2.5,0.5)(3,1.5)
\rput(1,0){${\scriptstyle \Delta_{irr}}$}
\end{pspicture}
\qquad\qquad
\begin{pspicture}(0,0)(2,3)
\pscurve(1.5,3)(1.25,2)(1,1.5)(0.5,1.2)
\pscurve(1.5,0)(1.318,0.75)(1.25,1)(1,1.5)(0.5,1.8)
\rput(1.45,2.06){${\scriptstyle I}$}
\rput(1.2,2.55){${\scriptstyle i}$}
\rput(1.45,0.99){${\scriptstyle J}$}
\rput(1.2,0.45){${\scriptstyle j}$}
\rput(1,-0.2){${\scriptstyle \Delta_{\{(i,I),(j,J)\}}}$}
\psdots(1.25,2)(1.31,2.2)(1.25,1)
\end{pspicture}
\]
Let $\delta_{irr}$ and $\delta_{\{(i,I),(j,J)\}}$ be the classes of
$\Delta_{irr}$ and $\Delta_{\{(i,I),(j,J)\}}$ in 
$\Pic(\bMM_{g,T}) \otimes \QQ$
respectively.
Strictly speaking, $\delta_{irr} = 
c_1(\OO_{\bMM_{g,T}}(\Delta_{irr}))$ and
\[
\delta_{\upsilon} =
\begin{cases}
{\displaystyle
\frac{1}{2} c_1\left(\OO_{\bMM_{g,T}}
\left(\Delta_{\upsilon}\right)\right)} & 
\text{if $\upsilon = \{ (1,\emptyset), (g-1, T) \}$}, \\
{ } \\
c_1\left(\OO_{\bMM_{g,T}}\left(\Delta_{\upsilon}\right)\right) & 
\text{otherwise}.
\end{cases}
\]
In the case where $T = \emptyset$, 
we denote $\delta_{\{(i,\emptyset),(j,\emptyset)\}}$
by $\delta_{\{i,j\}}$ or $\delta_{\min\{i,j\}}$, i.e.,
\[
\delta_i = \delta_{\{i,g-i\}} = \delta_{\{(i,\emptyset),(g-i,\emptyset)\}}
\qquad(i=1,\ldots,[g/2])
\]
on $\bMM_g$.

Let $\left(Z; \{ P_t \}_{t \in T} \right)$ 
be a $T$-pointed stable curve of genus 
$g$ over $k$.
Let $Q$ be a node of $Z$, and
$Z_Q$ the partial normalization of $Z$ at $Q$.
Then, the type of $Q$ is defined as follows:
\begin{enumerate}
\renewcommand{\labelenumi}{$\bullet$}
\item
The case where $Z_Q$ is connected. Then, $Q$ is of type $0$.

\item
The case where $Z_Q$ is not connected.
Let $Z_1$ and $Z_2$ be two connected components of $Z_Q$.
Let $i$ (resp. $j$) be the arithmetic genus of $Z_1$ (resp. $Z_2$).
Let $I = \{ t \in T \mid P_t  \in Z_1 \}$  and 
$J = \{ t \in T \mid P_t \in Z_2 \}$.
Then, we say $Q$ is of type $\{(i,I), (j,J)\}$.
\end{enumerate}
In the case where $T = \emptyset$,
for simplicity,
a node of type $\{(i,\emptyset),(j,\emptyset)\}$
is said to be of type $i$, where $i \leq j$.

Let $Y$ be a normal variety, 
and let $f : X \to Y$ be a $T$-pointed
stable curve of genus $g$ over $Y$.
Let $Y_0$ be the maximal open set over which $f$ is smooth.
Assume that $Y_0 \not= \emptyset$.
For $x \in X$, we define $\mult_x(X)$ to be 
$\length_{\OO_{X,x}}(\omega_{X/Y}/\Omega_{X/Y})$.
If $x$ is the generic point of a subvariety $T$, then
we denote $\mult_x(X)$ by $\mult_T(X)$.
If $x$ is closed, $Y$ is smooth at $f(x)$ and $Y \setminus Y_0$ 
is smooth at $f(x)$,
then $X$ is locally given by $\{ xy = t^{\mult_x(X)} \}$ around $x$,
where $t$ is a defining equation of $Y \setminus Y_0$ around $f(x)$.
Thus, if $Y$ is a curve, then the type of  singularity at $x$ is 
$A_{\mult_x(X)-1}$.

Here, for $\upsilon \in  \overline{\Upsilon}_{g,T}$,
let $S(X/Y)_{\upsilon}$ (resp.
$S(X/Y)_{irr}$)  be the set of irreducible components of $\Sing(f)$ such that
the type of $s$ in $f^{-1}(f(s))$ for a general 
$s \in S(X/Y)_{\upsilon}$ (resp.
$S(X/Y)_{irr}$) is $\upsilon$
(resp. $0$).
We set
\[
\delta_{\upsilon}(X/Y) =
\sum_{S \in S(X/Y)_{\upsilon}} \mult_{S}(X) f_*(S)
\]
and
\[
\delta_{irr}(X/Y) =
\sum_{S \in S(X/Y)_{irr}} \mult_{S}(X) f_*(S).
\]
Then, $\delta_{irr}$ and $\delta_{\upsilon}$ 
are normalized to guarantee the following formula:
\[
\delta_{irr}(X/Y) = \varphi^*(\delta_{irr})
\quad\text{and}\quad
\delta_{\upsilon}(X/Y) = 
\varphi^*(\delta_{\upsilon})
\]
in $A^1(Y) \otimes \QQ$,
where $\varphi : Y \to \bMM_{g,T}$ is the induced morphism by
$X \to Y$.

\subsection{The clutching maps}
\setcounter{Theorem}{0}
\label{subsec:clutching:maps}
Here let us consider the clutching maps and their properties.

Let $\pi : X \to Y$ be a prestable curve, i.e.,
$\pi : X \to Y$ is a flat and proper morphism such that
the geometric fibers of $\pi$ are reduced curves with
at most ordinary double points. We don't assume the connectedness of fibers.
Let $s_1, s_2 : Y \to X$ be two non-crossing sections such that
$\pi$ is smooth at points $s_1(y)$ and $s_2(y)$ ($\forall y \in Y$).
Then, by virtue of \cite[Theorem~3.4]{K23},
we have the clutching diagram:
\[
\xymatrix{
X \ar[rr]^{p} \ar[dr]^{\pi} & & X' \ar[dl]_{\pi'} \\
& Y \ar@/^0.7pc/[ul]|{s_1} \ar@/^1.4pc/[ul]|{s_2} \ar@/_0.7pc/[ur]|{s} &
}
\]
Roughly speaking, $X'$ is a prestable curve over $Y$ obtained by
identifying $s_1(Y)$ with $s_2(Y)$, and
$s$ is a section of $X' \to Y$ with $p \cdot s_1 = p \cdot s_2 = s$.
For details, see \cite[Theorem~3.4]{K23}.

We assume that $\pi' : X' \to Y$ is a $T$-pointed stable curve of genus $g$,
and $s$ is one of sections of $\pi' : X' \to Y$ arising from
$T$-points of $\pi' : X' \to Y$.
Let $\varphi : Y \to \bMM_{g,T}$ be the induced morphism.
Here we set $\Lambda = \det(R\pi_*(\omega_{X/Y}))$,
$\Delta = \det(R\pi_*(\omega_{X/Y}/\Omega_{X/Y}))$ and
$\Psi = s_1^*(\omega_{X/Y}) \otimes s_2^*(\omega_{X/Y})$.
Then, we have the following.

\begin{Proposition}
\label{prop:clutch:map:line:bundle}
For simplicity, the divisor $\delta_{irr}$ on $\bMM_{g,T}$
is denoted by $\delta_0$.
\begin{enumerate}
\renewcommand{\labelenumi}{(\arabic{enumi})}
\item
$\varphi^{*}(\lambda) = \Lambda$ and $\varphi^*(\delta) = -\Psi + \Delta$,
where $\delta = \delta_0 + \sum_{\upsilon \in 
\overline{\Upsilon}_{g,T}} 
\delta_{\upsilon}$.

\item
We assume that $\pi(\Sing(\pi)) \not= Y$ and
every geometric fiber of $\pi$ has one node at most.
Let
\[
\Delta = \Delta_0 + \sum_{\upsilon \in \overline{\Upsilon}_{g,T}} 
\Delta_{\upsilon}
\] 
be the decomposition such that
the node of $\pi^{-1}(x)$ \rom{(}$x \in (\Delta_{t})_{red}$\rom{)} 
gives rise to
a node of type $t$ in ${\pi'}^{-1}(x)$.
Moreover, let $a$ be the type of 
$s(y)$ in ${\pi'}^{-1}(y)$ \rom{(}$y \in Y$\rom{)}.
Then,
\[
\varphi^{*}(\delta_t) = \begin{cases}
-\Psi + \Delta_a & \text{if $t = a$} \\
\Delta_t & \text{if $t \not= a$}
\end{cases}
\]
\end{enumerate}
\end{Proposition}

\Proof
(1)
Since $\det(R{\pi'}_*(\omega_{X'/Y})) = \det(R\pi_*(\omega_{X/Y}))$, 
the first statement is obvious.
Thus, we can see that
\begin{align*}
\varphi^*(\delta) & = \det(R{\pi'}_*(\omega_{X'/Y}/\Omega_{X'/Y})) \\
& = \det(R{\pi'}_*(\omega_{X'/Y})) -
\det(R{\pi'}_*(\Omega_{X'/Y})) \\
& = \Lambda - \det(R{\pi'}_*(\Omega_{X'/Y})).
\end{align*}
On the other hand, by \cite[Theorem~3.5]{K23},
there is an exact sequence
\[
0 \to s^*(\Psi) \to \Omega_{X'/Y} \to p_*(\Omega_{X/Y}) \to 0.
\]
Therefore, we get (1).

\medskip
(2) This is a consequence of (1).
\QED

As a corollary, we have the following.

\begin{Corollary}
\label{cor:formula:clutching:map}
\begin{enumerate}
\renewcommand{\labelenumi}{(\arabic{enumi})}
\item
Let $a$ and $b$ be non-negative integers, and
$T$ and $S$ non-empty finite sets with $T \cap S = \emptyset$,
$2a -2 + \vert T \vert > 0$ and $2b - 2 + \vert S \vert > 0$.
Let us fix $s \in S$ and $t \in T$, and set $T' = T \setminus \{t\}$ and
$S' = S \setminus \{ s \}$.
Let $\alpha : \bMM_{a,T} \times \bMM_{b,B}\to 
\bMM_{a+b, T' \cup S'}$ 
be the clutching map, and
$p : \bMM_{a,T} \times \bMM_{b,S} \to \bMM_{a,T}$ and
$q : \bMM_{a,T} \times \bMM_{b,S} \to \bMM_{b,S}$ 
the projection to the first factor and
the projection to the second factor respectively. 
We set divisors $D \in \Pic(\bMM_{a+b,T' \cup S'}) \otimes \QQ$,
$E \in \Pic(\bMM_{a,T}) \otimes \QQ$ and 
$F \in \Pic(\bMM_{b,S})\otimes \QQ$ as follows:
\begin{align*}
D & = c \lambda + \sum_{l \in T' \cup S'} d_l \psi_l 
+ e_{irr} \delta_{irr} + \sum_{\{ (i, I), (j, J) \} \in 
\overline{\Upsilon}_{a+b, T' \cup S'}} e_{\{(i, I),(j,J)\}}
\delta_{\{(i,I),(j,J)\}}, \\
E & =c\lambda - e_{\{(a, T'), (b, S')\}}
\psi_{t} + \sum_{l \in T'} d_l \psi_l 
+ e_{irr} \delta_{irr} \\
& \phantom{c\lambda - e_{\{(a, T'), (b, S')\}}\psi_{t}} 
+ \sum_{\{ (i', I'), (j', J') \} \in 
\overline{\Upsilon}_{a, T} \atop t \in J'} 
e_{\{ (i', I'), (j'+b, J' \cup S'\setminus \{t \}) \}}
\delta_{\{ (i', I'), (j', J') \}}, \\
F & =
c\lambda - e_{\{(a, T'), (b, S')\}}
\psi_{s} + \sum_{l \in S'} d_l \psi_l 
+ e_{irr} \delta_{irr} \\
& \phantom{c\lambda - e_{\{(a, T'), (b, S')\}}\psi_{s}}
 + \sum_{\{ (i'', I''), (j'', J'') \} \in 
\overline{\Upsilon}_{b, S} \atop s \in J''} 
e_{\{ (i'', I''), (j''+a, J'' \cup T' \setminus \{s\}) \}}
\delta_{\{ (i'', I''), (j'', J'') \}}.
\end{align*}
Then $\alpha^*(D) = p^*(E) + q^*(F)$.

\item
Let $g$ be a non-negative integer and $T$ a finite set with
$\vert T \vert \geq 2$ and $2g -2 + \vert T \vert > 0$.
Let us fix two elements $t, t' \in T$, and set $T' = T \setminus \{t,t'\}$.
Let $\beta : \bMM_{g,T} \to \bMM_{g+1, T'}$ be the clutching map.
We set $D \in \Pic(\bMM_{g+1,T'}) \otimes \QQ$ and $E \in 
\Pic(\bMM_{g,T}) \otimes \QQ$
as follows:
\begin{align*}
D & = c \lambda + \sum_{l \in T'} d_l \psi_l 
+ e_{irr} \delta_{irr} + \sum_{\{ (i, I), (j, J) \} \in 
\overline{\Upsilon}_{g+1, T'}} e_{\{(i, I),(j,J)\}}
\delta_{\{(i,I),(j,J)\}}, \\
E & =c\lambda - e_{irr}(\psi_t + \psi_{t'})
+ \sum_{l \in T'} d_l \psi_l 
+ e_{irr} \delta_{irr} + 
\sum_{\{ (i', I'), (j', J') \} \in 
\overline{\Upsilon}_{g, T} \atop t \in I', t' \in J'} 
e_{irr} \delta_{\{ (i', I'), (j', J') \}} \\
& \phantom{c\lambda - e_{irr}(\psi_t + \psi_{t'})} 
+ \sum_{\{ (i', I'), (j', J') \} \in 
\overline{\Upsilon}_{g, T} \atop t, t' \in J'} 
e_{\{ (i', I'), (j'+1, J' \setminus \{t,t' \}) \}}
\delta_{\{ (i', I'), (j', J') \}}.
\end{align*}
Then $\beta^*(D) = E$.
\end{enumerate}
\end{Corollary}

\Proof
In the following, for $x \in \bMM_{*,*}$,
we denote by $C_x$ the corresponding nodal curve to $x$.

\medskip
(1) If $C_{\alpha(x,y)}$ has two nodes, then
we denote by $ty(x,y)$ the type of the node different from
the node arising from the clutching map.
Then, 
\begin{multline*}
ty(x,y) = \\
\begin{cases}
\{ (i',I'), (j'+b, J' \cup S' \setminus\{ t \})\} & 
\text{if $x \in \Delta_{\{(i',I'),(j',J')\}} \cap \bMM_{a,T}^{[1]}$, 
$y \in \MM_{b,S}$ and $t \in J'$}, \\
\{ (i'',I''), (j''+a, J'' \cup T' \setminus\{ s \})\} & 
\text{if $x \in \MM_{a,T}$,
$y \in \Delta_{\{(i'',I''),(j'',J'')\}} \cap \bMM_{b,S}^{[1]}$ and
$s \in J''$}.
\end{cases}
\end{multline*}
Thus, we get (1) by the above proposition.

\medskip
(2) In the same way as above,
if $C_{\beta(x)}$ has two nodes, then
we denote by $ty'(x)$ the type of the node different from
the node arising from the clutching map.
Then,
\[
ty'(x) = \begin{cases}
0 & \text{if $x \in \left(\Delta_{irr} \cup 
\bigcup_{t \in I', t' \in J'} 
\Delta_{\{ (i',I'), (j',J') \}} \right) \cap \bMM_{g,T}^{[1]}$}, \\
\{ (i', I'), (j'+1, J'\setminus\{t,t'\}) \} & 
\text{if $x \in \Delta_{\{(i',I'),(j',J')\}} \cap \bMM_{g,T}^{[1]}$
and $t,t' \in J'$},
\end{cases}
\]
which implies (2) by the above proposition.
\QED

\renewcommand{\theTheorem}{\arabic{section}.\arabic{Theorem}}
\renewcommand{\theClaim}{\arabic{section}.\arabic{Theorem}.\arabic{Claim}}
\renewcommand{\theequation}{\arabic{section}.\arabic{Theorem}.\arabic{Claim}}

\section{A generalization of relative Bogomolov's inequality}
Let $f : X \to Y$ be a projective morphism of quasi-projective varieties
of the relative dimension one, and let $E$ be a locally free sheaf on
$X$. Let us fix a point $y \in Y$.
Assume that $f$ is smooth over $y$ and $\rest{E}{f^{-1}(y)}$
is strongly semistable. In the paper \cite{MoRB}, we proved
that $\dis_{X/Y}(E)$ is weakly positive at $y$
under the assumption that $Y$ is smooth.
In this section, we generalize it to the case where $Y$ is normal.

\begin{Proposition}[$\ch(k) \geq 0$]
\label{prop:RR:formula}
Let $X$ and $Y$ be algebraic varieties, and
$f : X \to Y$ a surjective and projective morphism of $\dim f = d$.
Let $L$ and $A$ be line bundles on $X$. If $Y$ is normal, then
there are $\QQ$-divisors $Z_0, \ldots, Z_d$ on $Y$ such that
\[
c_1\left( Rf_*(L^{\otimes n} \otimes A) \right) \sim_{\QQ}
\frac{ f_*(c_1(L)^{d+1})}{(d+1)!} n^{d+1} + 
\sum_{i=0}^{d} Z_i n^i
\]
for all $n > 0$.
\end{Proposition}

\Proof
We set $Y^0 = Y \setminus \Sing(Y)$, $X^0 = f^{-1}(Y^0)$ and
$f^0 = \rest{f}{X^0}$. Then, we have
\[
 c_1\left(Rf^0_*\left(\rest{(L^{\otimes n} \otimes A)}{X^0}\right)\right) =
 \rest{c_1(Rf_*(L^{\otimes n} \otimes A))}{Y^0}
\]
and
\[
 f^0_*\left(c_1\left(\rest{L}{X^0}\right)^{d+1}\right) =
 \rest{f_*(c_1(L)^{d+1})}{Y^0}.
\]
Thus, by virtue of \cite[Lemma~2.3]{MoRB},
there are $\QQ$-divisors $Z^0_0, \ldots, Z^0_d$ on $Y^0$ such that
\[
\rest{c_1(Rf_*(L^{\otimes n} \otimes A))}{Y^0} \sim_{\QQ}
\frac{\rest{f_*(c_1(L)^{d+1})}{Y^0}}{(d+1)!} n^{d+1} + 
\sum_{i=0}^{d} Z^0_i n^i
\]
for all $n > 0$.
Let $Z_i$ be the Zariski closure of $Z_i^0$ in $Y$.
Then, since $\codim(\Sing(Y)) \geq 2$,
\[
c_1(Rf_*(L^{\otimes n} \otimes A)) \sim_{\QQ}
\frac{f_*(c_1(L)^{d+1})}{(d+1)!} n^{d+1} + 
\sum_{i=0}^{d} Z_i n^i
\]
for all $n > 0$.
\QED

\begin{Theorem}[$\ch(k) \geq 0$]
\label{thm:nef:psudo:dis}
Let $X$ be a quasi-projective variety, $Y$ a normal quasi-projective
variety,
and $f : X \to Y$ a surjective and projective morphism of $\dim f = 1$.
Let $F$ be a locally free sheaf on $X$ with 
$f_*(c_1(F)) = 0$, and $S$ a finite subset of $Y$.
We assume that
$f$ is flat over any points of $S$, and that, for all $s \in S$, 
there are line bundles $L_{\bar{s}}$ and $M_{\bar{s}}$ on
the geometric fiber $X_{\bar{s}}$ over $s$
such that
\[
H^0(X_{\bar{s}}, \Sym^m(F_{\bar{s}}) \otimes L_{\bar{s}}) = 
H^1(X_{\bar{s}}, \Sym^m(F_{\bar{s}}) \otimes M_{\bar{s}}) = 0
\]
for $m \gg 0$.
Then, $f_*\left( c_2(F) - c_1(F)^2 \right)$ is weakly positive over $S$.
\end{Theorem}

\Proof
The proof of this theorem is exactly same as \cite[Theorem~2.4]{MoRB}
using Proposition~\ref{prop:RR:formula}, Proposition~\ref{prop:samp:det} and
\cite[Proposition~2.2]{MoRB}.
For reader's convenience, we give the sketch of the proof of it.

\medskip
Let $A$ be a very ample line bundle on $X$ such that 
$A_{\bar{s}} \otimes L_{\bar{s}}$ and $A_{\bar{s}} 
\otimes M_{\bar{s}}^{\otimes -1}$ 
are very ample on $X_{\bar{s}}$ for all $s \in S$. 
Then, we can see the following claim in the same way as in
\cite[Claim~2.4.1]{MoRB}

\begin{Claim}
\label{claim:thm:nef:psudo:dis:0}
$H^0(X_s, \Sym^m(F_s) \otimes A_s^{\otimes -1}) = 
H^1(X_s, \Sym^m(F_s) \otimes A_s) = 0$
for all $s \in S$ and $m \gg 0$.
\end{Claim}

Since $X$ is an integral scheme 
of dimension greater than or equal to $2$,
and $X_s$ ($s \in S$) is a $1$-dimensional scheme over $\kappa(s)$,
there is $B \in |A^{\otimes 2}|$ such that
$B$ is integral, and that $B \cap X_s$ is finite for all $s \in S$, i.e.,
$B$ is finite over any points of $S$.
Let $\pi : B \to Y$ be the morphism induced by $f$.
Let $H$ be an ample line bundle on $Y$ such that
$\pi_*(F_B)  \otimes H$ and $\pi_*(A_B) \otimes H$ are generated 
by global sections at
any points of $S$,
where $F_B = \rest{F}{B}$ and $A_B = \rest{A}{B}$.

By using Proposition~\ref{prop:RR:formula},
there are $\QQ$-divisors $Z_0, \ldots, Z_r$ on $Y$
such that
\begin{multline*}
\sum_{i \geq 0}
(-1)^i c_1\left( R^if_*\left(\Sym^{m}(F \otimes f^*(H)) \otimes 
A^{\otimes -1} \otimes f^*(H) 
\right) \right) \\
\sim_{\QQ} -\frac{1}{(r+1)!} f_*  (c_2(F) - c_1(F)^2)  m^{r + 1}
+ \sum_{i=0}^{r} Z_i m^i
\end{multline*}
in the same way as in the proof of \cite[Theorem~2.4]{MoRB}.
The following claim also can be proved in the same way as in
 \cite[Claim~2.4.2]{MoRB}.

\begin{Claim}
\label{claim:thm:nef:psudo:dis:1}
If $m \gg 0$, then we have the following.
\begin{enumerate}
\renewcommand{\labelenumi}{(\alph{enumi})}
\item
$c_1\left( R^if_*\left(\Sym^{m}(F \otimes f^*(H)) \otimes 
A^{\otimes -1} \otimes f^*(H)
\right) \right) = 0$
for all $i \geq 2$.

\item
$f_*\left(\Sym^{m}(F \otimes f^*(H)) \otimes A^{\otimes -1} 
\otimes f^*(H) \right) = 0$.

\item
$R^1f_*\left(\Sym^{m}(F \otimes f^*(H)) \otimes A^{\otimes -1} 
\otimes f^*(H) \right)$
is free at any points of $S$.

\item
$R^1 f_*\left(\Sym^{m}(F \otimes f^*(H)) \otimes A \otimes f^*(H) \right) = 0$ 
around any points of $S$.
\end{enumerate}
\end{Claim}

By (a) and (b) of Claim~\ref{claim:thm:nef:psudo:dis:1},
\[
\frac{f_* (c_2(F) - c_1(F)^2)}{(r+1)!} 
\sim_{\QQ} \frac{c_1\left( R^1 f_*\left(\Sym^{m}(F \otimes f^*(H)) 
\otimes A^{\otimes -1} \otimes f^*(H)
\right) \right)}{m^{r+1}}
+ \sum_{i=0}^{r} \frac{Z_i}{m^{r+1-i}}.
\]
Hence, it is sufficient to show that 
\[
c_1\left( R^1 f_*\left(\Sym^{m}(F \otimes f^*(H)) \otimes A^{\otimes -1} 
\otimes f^*(H)
\right) \right)
\]
is semi-ample over $S$.
This can be proved in the same way as in the proof of \cite[Theorem~2.4]{MoRB}
by using \cite[Proposition~2.2]{MoRB},
Claim~\ref{claim:thm:nef:psudo:dis:1}
and Proposition~\ref{prop:samp:det}.
\QED

Let $C$ be a smooth projective curve and $E$ a vector bundle on $C$.
We say $E$ is {\em strongly semistable} if,
for any finite morphisms $\phi : C' \to C$ of
smooth projective curves, $\phi^*(E)$ is semistable.
Note that if $\ch(k) = 0$ and $E$ is semistable, then
$E$ is strongly semistable.
As a corollary,
we have the following, which can be proved in the exactly same way
as \cite[Corollary~2.5]{MoRB}.

\begin{Corollary}[$\ch(k) \geq 0$]
\label{cor:rel:bogo}
Let $X$ be a quasi-projective variety, 
$Y$ a normal quasi-projective variety,
and $f : X \to Y$ a surjective and projective morphism 
of $\dim f = 1$.
Let $E$ be a locally free sheaf on $X$ and
$S$ a finite subset of $Y$.
If, for all $s \in S$, $f$ is flat over $s$, 
the geometric fiber $X_{\bar{s}}$ over $s$ is reduced and Gorenstein,
and $E$ is strongly semistable on each
connected component of the normalization of
$X_{\bar{s}}$, 
then $\dis_{X/Y}(E)$ is weakly positive over $S$.
\end{Corollary}

\begin{Remark}[$\ch(k) = 0$]
\label{rem:slope:inq:mg}
In \cite{MoRB}, we proved that the divisor
\[
 (8g+4)\lambda - g\delta_{irr} - \sum_{i=1}^{[g/2]}4i(g-i)\delta_i
\]
on $\bMM_g$ is weakly positive over any finite subsets of $\MM_g$.
Here we give an alternative proof of this inequality.

Fix a polynomial $P_g(m) = (6m-1)(g-1)$.
Let $H_g \subset \Hilb^{P_g}_{\PP^{5g-6}}$ be a subscheme
of all tri-canonically embedded stable curves, 
$Z_g \subset H_g \times \PP^{5g-6}$ 
the universal tri-canonically embedded stable curves,
and $f_g : Z_g \to H_g$ the natural projection.
Then, $G = \PGL(5g-5)$ acts on $Z_g$ and $H_g$, and
$f_g$ is a $G$-morphism.
Let $\phi : H_g \to \bMM_g$ be the natural morphism
of the geometric quotient.
Then, by Seshadri's theorem \cite[Theorem~6.1]{Se},
there is a finite morphism $h : Y \to \bMM_g$ of normal varieties
with the following properties.
Let $W_g$ be the normalization of $H_g \times_{\bMM_g} Y$, and
let $\pi : W_g \to H_g$ and $\phi' : W_g \to Y$ be the induced morphisms
by the projections of 
$H_g \times_{\bMM_g} Y \to H_g$ and $H_g \times_{\bMM_g} Y \to Y$
respectively. Then, we have the following.
\begin{enumerate}
\renewcommand{\labelenumi}{(\arabic{enumi})}
\item
$G$ acts on $W_g$, and $\pi$ is a $G$-morphism.

\item
$\phi' : W_g \to Y$ is a principal $G$-bundle.
\end{enumerate}
Thus, $f'_g : U_g = Z_g \times_{H_g} W_g \to W_g$ is a stable curve,
$G$ acts on $U_g$ and $f'_g$ is a $G$-morphism.
Since $\phi' : W_g \to Y$ is a principal $G$-bundle,
we can easily see that $U_g$ is also a principal $G$-bundle and the geometric
quotient $X = U_g/G$ gives rise to a stable curve
$f : X \to Y$ over $Y$. Moreover, $U_g = W_g \times_Y X$.
Then, we have the following commutative diagram:
\[
\xymatrix@-0.8pc{
& Z_g \ar[ld]_{f_g} & & U_g \ar[ll]_{\pi'} \ar[ld]_{f'_g} \ar[dd]^{\phi''} \\
H_g \ar[dd]_{\phi} & & W_g \ar[ll]^{\pi} \ar[dd]_{\phi'} & \\
& & & X \ar[ld]^{f} \\
\bMM_g & & Y \ar[ll]^{h} &
}
\]

Let $\Delta$ be the minimal closed subset of $H_g$ such that
$f_g$ is not smooth over a point of $\Delta$.
Then, by \cite[Theorem~(1.6) and Corollary~(1.9)]{DM},
$Z_g$ and $H_g$ are quasi-projective and smooth, 
and $\Delta$ is a divisor with only normal crossings.
Let $\Delta = \Delta_{irr} \cup \Delta_1 \cup \cdots \cup \Delta_{[g/2]}$
be the irreducible decomposition of $\Delta$ such that,
if $x \in \Delta_i \setminus \Sing(\Delta)$
(resp. $x \in \Delta_{irr} \setminus \Sing(\Delta)$), then
$f_g^{-1}(x)$ is a stable curve with one node of type $i$
(resp. irreducible stable curve with one node).

Form now on, we consider everything over $\bMM_g^{[1]}$.
(Recall that $\bMM_g^{[1]}$ is the set of stable curves with one node at most.)
In the following, the superscript ``0'' means the objects over $\bMM_g^{[1]}$.

In \cite[\S3]{Mo5}, 
we constructed a locally sheaf $F$ on $Z^0_g$ with the following
properties.
\begin{enumerate}
\renewcommand{\labelenumi}{(\alph{enumi})}
\item
$F$ is invariant by the action of $G$.

\item
For each $y \in H^0_g \setminus (\Delta_1 \cup \cdots \cup \Delta_{[g/2]})$, 
$\rest{F}{f_g^{-1}(y)} = 
\Ker\left(H^0(\omega_{f_g^{-1}(y)}) \otimes \OO_{f_g^{-1}(y)} \to
\omega_{f_g^{-1}(y)}\right)$, which is semistable on $f_g^{-1}(y)$.

\item
$\dis_{Z^0_g/H^0_g}(F) = (8g+4) \det(\pi_*(\omega_{Z^0_g/H^0_g})) - 
g \Delta^0_{irr} - \sum_{i=1}^{\left[\frac{g}{2}\right]} 4i(g-i) \Delta^0_i$.
\end{enumerate}

Then, ${\pi'}^*(F)$ is a $G$-invariant locally free sheaf on $U_g^0$, so that
${\pi'}^*(F)$ can be descended to $X^0$ because
$U_g \to X$ is a principal $G$-bundle. Namely,
there is a locally free sheaf $F'$ on $X^0$ such that
${\phi''}^*(F') = {\pi'}^*(F)$.
Therefore, by Lemma~\ref{lem:comp:dis},
${\phi'}^*(\dis_{X^0/Y^0}(F')) = \pi^*(\dis_{Z_g^0/H_g^0}(F))$.
On the other hand, if we set
\[
D = (8g+4) \lambda - 
g \delta_{irr} - \sum_{i=1}^{\left[\frac{g}{2}\right]} 4i(g-i) \delta_i,
\]
then $\phi^*(D^0) = \dis_{Z_g^0/H_g^0}(F)$.
Therefore, we get ${\phi'}^*(h^*(D^0)) = {\phi'}^*(\dis_{X^0/Y^0}(F'))$,
which implies that $h^*(D^0) = \dis_{X^0/Y^0}(F')$ because
$\Pic(W_g)^G = \Pic(Y)$.
Moreover, by Corollary~\ref{cor:rel:bogo},
$\dis_{X^0/Y^0}(F')$ is weakly positive over any finite subsets of
$h^{-1}(\MM_g)$.
Thus, $h_*(\dis_{X^0/Y^0}(F')) = \deg(h)D^0$ is
weakly positive over any finite subsets of $\MM_g$
by (2) of Proposition~\ref{prop:samp:pamp:push}.
Hence, $D$ is weakly positive over any finite subsets of of $\MM_g$
because $\codim(\bMM_g \setminus \bMM_g^{[1]}) \geq 2$.
\end{Remark}

\section{A certain kind of hyperelliptic fibrations}
We say $f : X \to Y$ is a hyperelliptic fibered surface
of genus $g$ if $X$ is a smooth projective surface,
$Y$ is a smooth projective curve,
the generic fiber of $f$ is a smooth hyperelliptic curve of genus $g$, and
there is no $(-1)$-curve along each fiber of $f$.
Since $f$ is minimal, the hyperelliptic involution of
the generic fiber extends to an automorphism
of $X$ over $Y$. We denote this automorphism by $j$.
Clearly, the order of $j$ is $2$, namely, $j \not= \operatorname{id}_X$
and $j^2 = \operatorname{id}_X$.

The purpose of this section is to show the existence of a special kind
of hyperelliptic fibered surfaces as described in the following propositions.

\begin{Proposition}[$\ch(k) = 0$]
\label{prop:a:ga:one}
For fixed integers $g$ and $i$ with $g \geq 2$ and $0 \leq i \leq g-1$,
there is a hyperelliptic fibered surface $f : X \to Y$ of genus $g$, and
a section $\Gamma$ of $f$ such that
\begin{enumerate}
\renewcommand{\labelenumi}{(\arabic{enumi})}
\item
$\Sing(f) \not= \emptyset$, $j(\Gamma) = \Gamma$,

\item
every singular fiber of $f$ is a stable curve consisting 
a smooth projective curve of genus $i$
and a smooth projective curve of genus $g-i$ 
meeting transversally at one point, and that

\item
$\Gamma$ intersects with the component of genus $g-i$ on every singular fiber.
\end{enumerate}
\end{Proposition}

\begin{Proposition}[$\ch(k) = 0$]
\label{prop:a:ga:two}
For fixed integers $g$ and $i$ with $g \geq 2$ and $0 \leq i \leq g$,
there is a hyperelliptic fibered surface $f : X \to Y$ of genus $g$, and
a section $\Gamma$ of $f$ such that
\begin{enumerate}
\renewcommand{\labelenumi}{(\arabic{enumi})}
\item
$\Sing(f) \not= \emptyset$,
$j(\Gamma) \cap \Gamma = \emptyset$,

\item
every singular fiber of $f$ is a stable curve 
consisting a smooth projective curve of genus $i$
and a smooth projective curve of genus $g-i$ 
meeting transversally at one point, and that

\item
$\Gamma$ intersects with the component of genus $g-i$ on every singular fiber.
\end{enumerate}
\end{Proposition}

\begin{Proposition}[$\ch(k) = 0$]
\label{prop:irr:one}
For fixed integers $g$ and $i$ with $g \geq 2$ and $0 \leq i \leq g-1$,
there is a hyperelliptic fibered surface $f : X \to Y$ of genus $g$, and
a section $\Gamma$ of $f$ such that
\begin{enumerate}
\renewcommand{\labelenumi}{(\arabic{enumi})}
\item
$\Sing(f) \not= \emptyset$, $j(\Gamma) = \Gamma$,

\item
every singular fiber of $f$ is a semi-stable curve consisting 
a smooth projective curve of genus $i$
and a smooth projective curve of genus $g-i-1$ 
meeting transversally at two points, and that

\item
$\Gamma$ intersects with the component of 
genus $g-i-1$ on every singular fiber.
\end{enumerate}
\end{Proposition}

\begin{Proposition}[$\ch(k) = 0$]
\label{prop:irr:two}
For fixed integers $g$ and $i$ with $g \geq 2$ and $0 \leq i \leq g-1$,
there is a hyperelliptic fibered surface $f : X \to Y$ of genus $g$, and
a section $\Gamma$ of $f$ such that
\begin{enumerate}
\renewcommand{\labelenumi}{(\arabic{enumi})}
\item
$\Sing(f) \not= \emptyset$, $j(\Gamma) \cap \Gamma = \emptyset$,

\item
every singular fiber of $f$ is a semi-stable curve consisting 
a smooth projective curve of genus $i$
and a smooth projective curve of genus $g-i-1$ 
meeting transversally at two points, and that

\item
$\Gamma$ intersects with the component of genus 
$g-i-1$ on every singular fiber.
\end{enumerate}
\end{Proposition}

\begin{Proposition}[$\ch(k) = 0$]
\label{prop:g:a:two:section}
For fixed integers $g$ and $i$ with $g \geq 2$ and $1 \leq i \leq g-1$,
there is a hyperelliptic fibered surface $f : X \to Y$ of genus $g$, and
non-crossing sections $\Gamma_1$ and $\Gamma_2$ of $f$ such that
\begin{enumerate}
\renewcommand{\labelenumi}{(\arabic{enumi})}
\item
$\Sing(f) \not= \emptyset$,
$j(\Gamma_1) = \Gamma_1$, $j(\Gamma_2) = \Gamma_2$,

\item
every singular fiber of $f$ is a stable curve consisting 
a smooth projective curve of genus $i$
and a smooth projective curve of genus $g-i$ 
meeting transversally at one point,

\item
$\Gamma_1$ and $\Gamma_2$ gives rise to a $2$-pointed
stable curve $(f : X \to Y, \Gamma_1, \Gamma_2)$, and that

\item
the type of $x$ in $f^{-1}(f(x))$ as $2$-pointed
stable curve is $\{ (i,\{1\}), (g-i,\{2\}) \}$ for all $x \in \Sing(f)$.
\end{enumerate}
\end{Proposition}

Let us begin with the following lemma.

\begin{Lemma}[$\ch(k) = 0$]
\label{lem:conic:fibration}
For non-negative integers $a_1$ and $a_2$,
there are a morphism $f_1 : X_1 \to Y_1$ 
of smooth projective varieties,
an effective divisor $D_1$ on $X_1$, a line bundle $L_1$ on $X_1$,
a line bundle $M_1$ on $Y_1$, and
non-crossing sections $\Gamma_1$ and $\Gamma_2$ of $f_1 : X_1 \to Y_1$ 
with the following properties.
\begin{enumerate}
\renewcommand{\labelenumi}{(\arabic{enumi})}
\item
$\dim X_1 = 2$ and $\dim Y_1 = 1$.

\item
Let $\Sigma_1$ be the set of all critical values of $f_1$, i.e.,
$P \in \Sigma_1$ if and only if $f_1^{-1}(P)$ is a singular variety.
Then, for any $P \in Y_1 \setminus \Sigma_1$, $f_1^{-1}(P)$ is a smooth
rational curve.

\item
$\Sigma_1 \not= \emptyset$, and for any $P \in \Sigma_1$, 
$f_1^{-1}(P)$ is a reduced curve
consisting of two smooth rational curves $E_P^{1}$ and $E_P^{2}$ joined
at one point transversally.

\item
$D_1$ is smooth and
$\rest{f_1}{D_1} : D_1 \to Y_1$ is etale.

\item
$(E_P^{1} \cdot D_1) = a_1 + 1$ and $(E_P^{2} \cdot D_1) = a_2 + 1$ 
for any $P \in \Sigma_1$.
Moreover, $D_1$ does not pass through $E_P^1 \cap E_P^2$.

\item
There is a map $\sigma : \Sigma_1 \to \{ 1, 2 \}$ such that
\[
D_1 \in \left|
L_1^{\otimes a_1 + a_2 + 2} \otimes f_1^*(M_1) \otimes 
\OO_{X_1}\left(-\sum_{P \in \Sigma_1} (a_{\sigma(P)} + 1) 
E_P^{\sigma(P)}\right) \right|.
\]

\item
$\deg(M_1)$ is divisible by $(a_1 + 1)(a_2 + 1)$.

\item
\[
 \Gamma_1 \in \left| L_1 \otimes \OO_{X_1}
\left( -\sum_{P \in \Sigma_1 \atop \sigma(P) = 1}E_P^1 \right) \right|
\quad\text{and}\quad
\Gamma_2 \in \left| L_1 \otimes \OO_{X_1}
\left( -\sum_{P \in \Sigma_1 \atop \sigma(P) = 2}E_P^2 \right) \right|.
\]
Moreover,
\[
(D_1 \cdot \Gamma_1) = (D_1 \cdot \Gamma_2) = 0
\quad\text{and}\quad
 (E_P^i \cdot \Gamma_j) = \begin{cases}
0 & \text{if $i \not= j$}, \\
1 & \text{if $i = j$}.
\end{cases}
\]
\end{enumerate}
\end{Lemma}

\Proof
We can prove this lemma in the exactly same way as in
\cite[Lemma~A.1]{MoRB} with a slight effort.
We use the notation in \cite[Lemma~A.1]{MoRB}.
Let $F_1$ and $F_2$ be curves in $\PP^1_{(X,Y)} \times \PP^1_{(S,T)}$
defined by $\{ X = 0 \}$ and $\{ X = Y \}$ respectively.
Note that $F_1 = p^{-1}((0:1))$, $F_2 = p^{-1}((1:1))$,
$D'' = p^{-1}((1:0))$, $(D' \cdot F_1) = (D' \cdot F_2) = 1$,
$D' \cap F_1 = \{ Q_1 \}$ and $D' \cap F_2 = \{ Q_2 \}$.
Then, since $u_1^*(F_1) = p_1^{-1}((0:1))$ and
$u_1^*(F_2) = p_1^{-1}((1:1))$,
in \cite[Claim~A.1.3]{MoRB}, we can see that
each tangent of $u_1^*(D)$ at $Q_{i,j}$ ($i=1,2$) is
different from $u_1^*(F_i)$.

Let $\overline{\Gamma}_i$ be the strict transform of
$u_1^*(F_i)$ by $\mu_1 : Z_1 \to \PP^1_{(X,Y)} \times Y$. Then,
\[
\overline{\Gamma}_i \in \left| 
\mu_1^*(p_1^*(\OO_{\PP^1}(1))) \otimes 
\OO_{Z_1}\left(-\sum_{j} E_{i,j}\right) \right|.
\]
Thus, if we set $\Gamma_i = (\nu_1)_*(\overline{\Gamma}_i)$, then
we get our lemma.
\QED

In the following proofs, we use
the notation in \cite[Proposition~A.2 and Proposition~A.3]{MoRB}.

\medskip
{\bf The proof of Proposition~\ref{prop:a:ga:one}:}\ 
We apply Lemma~\ref{lem:conic:fibration} to the case where
$a_1 = 2i$ and $a_2 = 2g - 2i - 1$.
We replace $D_2$ by $D_2 + \Gamma_2$ and $a_2$ by $a_2 + 1$.
Then, (4), (5) and (6) hold for the new $D_2$ and $a_2$.
Thus, we can construct $f : X \to Y$ in exactly same way as in 
\cite[Proposition~A.2]{MoRB}.
Since $u_2^*(\Gamma_2)$ is the ramification locus of $\mu_3$,
$\overline{\Gamma} = h^*(u_2^*(\Gamma_2))_{red}$ is a section of
$f_3$. Thus, if we set $\Gamma = \nu_3(\overline{\Gamma})$,
then we have our desired example.

\medskip
{\bf The proof of Proposition~\ref{prop:a:ga:two}:}\ 
Applying Lemma~\ref{lem:conic:fibration} to the case where
$a_1 = 2i$ and $a_2 = 2g - 2i$,
we can construct $f : X \to Y$ in exactly same way as in 
\cite[Proposition~A.2]{MoRB}.
Here let us consider $u_2^*(\Gamma_2)$. Then, $u_2^*(\Gamma_2)$
is a section of $f_2$ such that
$u_2^*(\Gamma_2) \cap (D_2 + B) = \emptyset$,
$(u_2^*(\Gamma_2) \cdot \overline{E}_Q^1) = 0$ and
$(u_2^*(\Gamma_2) \cdot \overline{E}_Q^2) = 1$ for all $Q \in \Sigma_2$.
Here we set $\Gamma' = \nu_3(\mu_3^*(u_2^*(\Gamma_2)))$.
Then, since $\mu_3^*(u_2^*(\Gamma_2))$ does not intersect with
the ramification locus of $\mu_3$, $\Gamma'$ is etale over $Y$.
Moreover, we can see
$(\Gamma' \cdot C_Q^1) = 0$ and
$(\Gamma' \cdot C_Q^2) = 2$ for all $Q \in \Sigma_2$.
If $\Gamma'$ is not irreducible, 
then we choose $\Gamma$ as one of irreducible component of $\Gamma'$.
If $\Gamma'$ is irreducible, then
we consider $X \times_Y \Gamma \to \Gamma$ and the natural section of
$X \times_Y \Gamma \to \Gamma$.
Then we get our desired example.

\medskip
{\bf The proof of Proposition~\ref{prop:irr:one}:}\ 
We apply Lemma~\ref{lem:conic:fibration} to the case where
$a_1 = 2i+1$ and $a_2 = 2g - 2i -2$.
We replace $D_2$ by $D_2 + \Gamma_2$ and $a_2$ by $a_2 + 1$.
Then, (4), (5) and (6) hold for the new $D_2$ and $a_2$.
Note that $\deg(M_1)$ is even.
Thus, we can get a double covering $\mu : X \to X_1$ 
in exactly same way as in 
\cite[Proposition~A.3]{MoRB}.
Let $f : X \to Y_1$ be the induced morphism, and
$\Gamma = \mu^*(\Gamma_2)_{red}$.
Then, we have our desired example.

\medskip
{\bf The proof of Proposition~\ref{prop:irr:two}:}\ 
Applying Lemma~\ref{lem:conic:fibration} to the case where
$a_1 = 2i+1$ and $a_2 = 2g - 2i -1$,
we can get a double covering $\mu : X \to X_1$ in exactly same way as in 
\cite[Proposition~A.3]{MoRB}.
Let $f : X \to Y_1$ be the induced morphism and
$\Gamma' = \mu^*(\Gamma_2)$.
Then, $\Gamma'$ is etale over $Y_1$.
If $\Gamma'$ is not irreducible, 
then we choose $\Gamma$ as one of irreducible component of $\Gamma'$.
If $\Gamma'$ is irreducible, then
we consider $X \times_{Y_1} \Gamma \to \Gamma$ and the natural section of
$X \times_{Y_1} \Gamma \to \Gamma$.
Then we get our desired example.

\medskip
{\bf The proof of Proposition~\ref{prop:g:a:two:section}:}
We apply Lemma~\ref{lem:conic:fibration} to the case where
$a_1 = 2i-1$ and $a_2 = 2g - 2i-1$.
We replace $D_1$ by $D_1 + \Gamma_1$, $D_2$ by $D_2 + \Gamma_2$,
$a_1$ by $a_1 + 1$, and $a_2$ by $a_2 + 1$.
Then, (4), (5) and (6) hold for the new $D_1$, $D_2$, $a_1$ and $a_2$.
Thus, we can construct $f : X \to Y$ in exactly same way as in 
\cite[Proposition~A.2]{MoRB}.
Since $u_2^*(\Gamma_1)$ and $u_2^*(\Gamma_2)$ are 
the ramification locus of $\mu_3$,
$\overline{\Gamma}_1 = h^*(u_2^*(\Gamma_1))_{red}$ and
$\overline{\Gamma}_2 = h^*(u_2^*(\Gamma_2))_{red}$ are
sections of
$f_3$. Thus, if we set $\Gamma_1 = \nu_3(\overline{\Gamma}_1)$
and $\Gamma_2 = \nu_3(\overline{\Gamma}_2)$,
then we have our desired example.

\bigskip
Finally, let us consider the following two lemmas, which
will be used in the later section.

\begin{Lemma}[$\ch(k) \geq 0$]
\label{lem:self:2:zero}
Let $X$ be a smooth projective surface and $Y$ a smooth projective curve.
Let $f : X \to Y$ be a surjective morphism with connected fibers, and let
$L$ be a line bundle on $X$.
If $\rest{L}{X_{\eta}}$ gives rise to a torsion element of
$\Pic(X_{\eta})$ on the generic fiber $X_{\eta}$ of $f$
and $\deg(\rest{L}{F}) = 0$ for every irreducible component $F$ of fibers, 
then we have $(L^2) = 0$.
\end{Lemma}

\Proof
Replacing $L$ by $L^{\otimes n}$ ($n \not= 0$), we may assume that
$\rest{L}{X_{\eta}} \simeq \OO_{X_{\eta}}$.
Thus, $f_*(L)$ is a line bundle on $Y$, and
the natural homomorphism $f^*f_*(L) \to L$ is injective. 
Hence, there is an effective
divisor $E$ on $X$ such that $f^*f_*(L) \otimes \OO_X(E) \simeq L$.
Since  $f^*f_*(L) \to L$ is surjective on the generic fiber,
$E$ is a vertical divisor.
Moreover, $(E \cdot F) = 0$ for every irreducible component $F$ of fibers.
Therefore, by Zariski's lemma, $(E^2) = 0$.
Hence, $(L^2) = (E^2) = 0$. 
\QED

\begin{Lemma}[$\ch(k) \geq 0$]
\label{lem:base:pt:free}
Let $C$ be a smooth projective curve of genus $g \geq 2$.
Let $\vartheta$ be a line bundle on $C$ with $\vartheta^{\otimes 2} = \omega_C$.
Let $\Delta$ be the diagonal of $C \times C$, and
$p : C \times C \to C$ and $q : C \times C \to C$ 
the projection to the first factor
and the projection to the second factor respectively. Then,
$p^*(\vartheta^{\otimes n}) \otimes q^*(\vartheta^{\otimes n}) 
\otimes \OO_{C \times C}((n-1)\Delta)$
is generated by global sections for all $n \geq 3$.
\end{Lemma}

\Proof
Since $p^*(\vartheta^{\otimes n}) \otimes q^*(\vartheta^{\otimes n})$ 
is generated by global sections,
the base locus of
$p^*(\vartheta^{\otimes n}) \otimes q^*(\vartheta^{\otimes n}) \otimes 
\OO_{C \times C}((n-1)\Delta)$
is contained in $\Delta$. Moreover,
\[
\rest{p^*(\vartheta^{\otimes n}) \otimes q^*(\vartheta^{\otimes n}) \otimes 
\OO_{C \times C}((n-1)\Delta)}{\Delta}
\simeq \omega_C.
\]
Thus, it is sufficient to see that
\[
H^0(p^*(\vartheta^{\otimes n}) \otimes q^*(\vartheta^{\otimes n}) \otimes 
\OO_{C \times C}((n-1)\Delta))
\to
H^0(\rest{p^*(\vartheta^{\otimes n}) \otimes q^*(\vartheta^{\otimes n}) 
\otimes \OO_{C \times C}((n-1)\Delta)}{\Delta})
\]
is surjective.

We define $L_{n,i}$ to be
\[
L_{n,i} = p^*(\vartheta^{\otimes n}) \otimes q^*(\vartheta^{\otimes n})
\otimes \OO_{C \times C}(i\Delta).
\]
Then, it suffices to check
$H^1(L_{n, n-2}) = 0$ for the above assertion.
By induction on $i$, we will see that $H^1(L_{n, i}) = 0$ for 
$0 \leq i \leq n-2$.

First of all, note that $H^1(\vartheta^{\otimes n}) = 0$ for $n \geq 3$.
Thus,
\[
H^1(p^*(\vartheta^{\otimes n}) \otimes q^*(\vartheta^{\otimes n})) =
H^1(p_*(p^*(\vartheta^{\otimes n}) \otimes q^*(\vartheta^{\otimes n}))) =
H^1(\vartheta^{\otimes n}) \otimes H^1(\vartheta^{\otimes n}) = 0.
\]
Moreover, let us consider the exact sequence
\[
0 \to L_{n,i-1} \to L_{n,i} \to \rest{L_{n,i}}{\Delta} \to 0.
\]
Here since $\rest{L_{n,i}}{\Delta} \simeq \omega_C^{\otimes n - i}$,
$H^1(\rest{L_{n,i}}{\Delta}) = 0$ if $i \leq n-2$. 
Thus, by the hypothesis of induction,
we can see $H^1(L_{n,i}) = 0$.
\QED

\section{Slope inequalities on $\bMM_{g,T}$}
Let $g$ be a non-negative integer and
$T$ a finite set 
with $2g-2 + \vert T \vert > 0$.
Recall that
\begin{align*}
\Upsilon_{g,T} & = \{ (i, I) \mid 
\text{$i \in \ZZ$, $0 \leq i \leq g$ and $I \subseteq T$} \} \setminus 
\left( \{ (0, \emptyset) \} \cup \{ (0, \{t\}) \}_{t \in T} \right), \\
\overline{\Upsilon}_{g,T} & = \{ \{(i,I), (j,J)\} \mid 
(i,I), (j,J) \in \Upsilon_{g,T}, i+j=g, I \cap J = \emptyset, I\cup J = T \}.
\end{align*}
For a subset $L$ of $T$, let us introduce a function
$\gamma_L : \Upsilon_{g,T} \times \Upsilon_{g,T} \to \ZZ$ given by
\[
\gamma_L((i,I),(j,J)) = \left( \det
\begin{pmatrix} i & \vert L \cap I \vert \\ j & \vert L \cap J \vert 
\end{pmatrix} + \vert L \cap I \vert\right)
\left( \det
\begin{pmatrix} i & \vert L \cap I \vert \\ j & \vert L \cap J \vert 
\end{pmatrix} - \vert L \cap J \vert \right).
\]
Note that $\gamma_L((i,I),(j,J)) = \gamma_L((j,J),(i,I))$, so that
$\gamma_L$ gives rise to a function on $\overline{\Upsilon}_{g,T}$.
Further, a $\QQ$-divisor $\theta_L$ on $\bMM_{g,T}$ is defined to be
\[
\theta_L = 4(g-1+\vert L \vert)(g-1) \sum_{t \in L} \psi_t 
- 12 \vert L \vert^2 \lambda
+ \vert L \vert^2  \delta_{irr} 
- \sum_{\upsilon \in \overline{\Upsilon}_{g,T}}
4 \gamma_L(\upsilon)
\delta_{\upsilon}.
\]
Then, we have the following.
 
\begin{Theorem}[$\ch(k) \geq 0$]
\label{thm:slope:inq:m:g:n}
For any subset $L$ of $T$,
the divisor $\theta_L$
is weakly positive over any finite subsets of $\MM_{g,T}$.
\end{Theorem}

\Proof
Clearly, we may assume $T = [n]$ for some non-negative integer $n$.
Let us take an $n$-pointed stable curve $f : X \to Y$
such that the induced morphism $h : Y \to \bMM_{g,[n]}$ is
a finite and surjective morphism of normal varieties.
Let $Y_0$ be the maximal Zariski open set of $Y$ over which
$f$ is smooth. 
Let $Y \setminus Y_0 = B_1 \cup \cdots \cup B_s$
be the irreducible decomposition 
of $Y \setminus Y_0$.
By using \cite[Lemma~3.2]{DeJ}, 
we can take a Zariski open set $Y_1$ with the following properties.
\begin{enumerate}
\renewcommand{\labelenumi}{(\arabic{enumi})}
\item
$\codim(Y \setminus Y_1) \geq 2$ and $Y_0 \subseteq Y_1$.

\item
$Y_1$ is smooth at any points of $Y_1 \cap (Y \setminus Y_0)$.

\item
$f : \Sing(f) \cap f^{-1}(Y_1) \to f(\Sing(f)) \cap Y_1$ is an isomorphism, 
so that for all $y \in Y_1$, the number of nodes of $f^{-1}(y)$ is
one at most.

\item
There is a projective birational morphism $\phi : Z_1 \to X_1 = f^{-1}(Y_1)$
such that if we set $f_1 = f \cdot \phi$, then $Z_1$ is smooth at any points
of $\Sing(f_1) \cap f_1^{-1}(Y \setminus Y_0)$ and 
$f_1 : Z_1 \to Y_1$ is an $n$-pointed semi-stable curve.
Moreover, $\phi$ is an isomorphism over $X_1 \setminus \Sing(f)$.

\item
For each $l = 1, \ldots, s$, there is a $t_l$ such that 
$\mult_x(X) = t_l + 1$ for all
$x \in \Sing(f)$ with $f(x) \in B_l \cap Y_1$.
\end{enumerate}
Let $K_0$ be a subset of $\{ 1, \ldots, s \}$ such that
$f^{-1}(x)$ is irreducible for all $x \in B_l \cap Y_1$, and
let $K_1 = \{ 1, \ldots, s \} \setminus K_0$.
For each $l \in K_1$, there is a $(g_l, I_l), (h_l, J_l) \in 
\Upsilon_{g,[n]}$
such that the type of $x$ is  $\{ (g_l, I_l), (h_l, J_l) \}$
for all $x \in \Sing(f)$ with $f(x) \in B_l \cap Y_1$.
From now on, by abuse of notation, we denote $B_l \cap Y_1$ by $B_l$.
For $l \in K_1$,
$f_1^{-1}(B_l)$ has two essential components $T_l^1$ and $T_l^2$, and
the components of $(-2)$-curves $E_1, \ldots, E_{t_l}$ such that
$T_l^1 \to B_l$ is an $I_l$-pointed smooth curve of genus $g_l$ and 
$T_l^2 \to B_l$
is a $J_l$-pointed smooth curve of genus $h_l$.
Moreover, the numbering of $E_1, \ldots, E_{t_l}$
is arranged as the following figure:
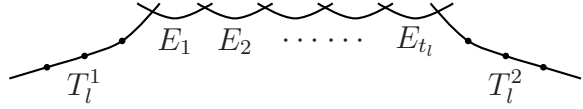
\[
\begin{pspicture}(0,0)(7,3)
\pscurve(0,1)(0.5,1.15)(1,1.3)(1.5,1.5)(2,2)
\pscurve(1.7,2)(2.2,1.8)(2.7,2)
\pscurve(2.5,2)(3,1.8)(3.5,2)
\pscurve(3.3,2)(3.8,1.8)(4.3,2)
\pscurve(4.1,2)(4.6,1.8)(5.1,2)
\pscurve(4.9,2)(5.4,1.8)(5.9,2)
\pscurve(7.6,1)(7.1,1.15)(6.6,1.3)(6.1,1.5)(5.6,2)
\psdots(6.1,1.5)(6.6,1.3)(7.1,1.15)
\psdots(0.5,1.15)(1,1.3)(1.5,1.5)
\rput(1,0.9){$T_l^1$}
\rput(6.6,0.9){$T_l^2$}
\rput(2.2,1.5){$E_1$}
\rput(3,1.5){$E_2$}
\rput(3.9,1.5){$\cdots$}
\rput(4.5,1.5){$\cdots$}
\rput(5.4,1.5){$E_{t_l}$}
\end{pspicture}
\]
Let $\Gamma_1, \ldots, \Gamma_n$ be the sections of
the $n$-pointed stable curve of $f : X \to Y$.
By abuse of notation, the lifting of $\Gamma_a$ to $Z_1$
is also denoted by $\Gamma_a$.
Here we consider a line bundle $L$ on $Z_1$ given by
\[
L = \omega_{Z_1/Y_1}^{\otimes \vert L \vert} \otimes 
\OO_{Z_1}\left(-(2g-2)\sum_{a \in L} \Gamma_a + 
\sum_{l \in K_1} \left( \vert L \vert (2g_l -1) - 
(2g-2) \vert L \cap I_l \vert \right)
\tilde{T}^1_l  \right),
\]
where
\[
\tilde{T}^1_l = (t_l + 1)T^1_l + \sum_{a=1}^{t_l} (t_l + 1 -a)E_a.
\]
We set $E = \OO_{X_1} \oplus L$.
Then, $\dis_{X_1/Y_1}(E) = -(f_1)_*(c_1(L)^2)$.
Here, we know the following formulae:
\[
\begin{cases}
f_*(c_1(\omega_{Z_1/Y_1}) \cdot \tilde{T}^1_l) = (t_l + 1)(2g_l-1)B_l, \\
f_*(\tilde{T}^1_l \cdot \tilde{T}^1_{l'}) = 
\begin{cases}
0 & \text{if $l\not=l'$}, \\
-(t_l+1)B_l & \text{if $l=l'$},
\end{cases} \\
f_*(\sum_{a \in L} \Gamma_a \cdot \tilde{T}^1_l) = 
(t_l+1)\vert L \cap I_l \vert B_l, \\
f_*(c_1(\omega_{Z_1/Y_1}) \cdot \Gamma_a) = 
-f_*(\Gamma_a \cdot \Gamma_a), \quad\text{(adjunction formula)} \\
12 \det(f_*(\omega_{Z_1/Y_1})) - \sum_{l=1}^s (t_l + 1)B_l = 
f_*(c_1(\omega_{Z_1/Y_1})^2)
\quad\text{(Noether's formula)}.
\end{cases}
\]
Thus, we can see that
\begin{multline*}
\dis_{Z_1/Y_1}(E) = 4(g-1+\vert L \vert)(g-1) 
f_*\left(c_1(\omega_{Z_1/Y_1}) \cdot 
\sum_{a \in L} \Gamma_a\right) - 
12\vert L \vert^2 \det(f_*(\omega_{Z_1/Y_1})) \\
+ \sum_{l \in K_0}\vert L \vert^2(t_l + 1)B_l
- \sum_{l \in K_1} 4(t_l + 1)\gamma_L(\{ (g_l,I_l),(h_l,J_l) \}) B_l.
\end{multline*}
On the other hand, 
for $y \in Y_0$, let $\phi : C' \to f^{-1}(y)$ be a finite morphism of
smooth projective curves. Then, $\phi^*(\rest{E}{f^{-1}(y)}) = 
\OO_{C'} \oplus \phi^*(\rest{L}{f^{-1}(y)})$ and
\[
\deg(\phi^*(\rest{L}{f^{-1}(y)})) = \deg(\phi) \deg(\rest{L}{f^{-1}(y)}) = 0.
\]
Therefore, $\phi^*(\rest{E}{f^{-1}(y)})$ is semistable,  which
means that $\rest{E}{f^{-1}(y)}$ is strongly semistable for all $y \in Y_0$.
Thus, by Corollary~\ref{cor:rel:bogo},
$\dis_{Z_1/Y_1}(E)$ is weakly positive over any finite subsets of
$Y_0$ as a divisor on $Y_1$. Therefore, if we set
\begin{multline*}
\theta'_L = 
4(g-1+\vert L \vert)(g-1) f_*\left(c_1(\omega_{Z_1/Y_1}) \cdot 
\sum_{a \in L} \Gamma_a\right) - 
12\vert L \vert^2 \det(f_*(\omega_{Z_1/Y_1})) \\
+ \sum_{l \in K_0}\vert L \vert^2(t_l + 1)B_l
- \sum_{l \in K_1} 4(t_l + 1)\gamma_L(\{ (g_l,I_l),(h_l,J_l) \}) B_l.
\end{multline*}
on $Y$, then $\theta'_L$ is weakly positive over any finite subsets of $Y_0$ 
as a divisor on $Y$.
Here 
$h^*(\theta_L) = \theta'_L$, so that $h_* (\theta'_L) = \deg(h) \theta_L$
by the projection formula.
Hence we have our theorem by (2) of
Proposition~\ref{prop:samp:pamp:push}.
\QED

Let us apply Theorem~\ref{thm:slope:inq:m:g:n}
to the cases $\bMM_{g,1}$ and $\bMM_{g,2}$.

\begin{Corollary}[$\ch(k) = 0$]
\label{cor:slope:cone:mg1}
Let $\bMM_{g,1} = \bMM_{g,\{1\}}$ be the moduli space of one-pointed
stable curves of genus $g \geq 1$.
We set $\delta_{i}, \mu, \theta_1 \in \Pic(\bMM_{g,1}) \otimes \QQ$
as follows:
\[
\begin{cases}
\delta_{i} = \delta_{\{(i, \emptyset), (g-i, \{1\})\}}
\quad(1 \leq i \leq g-1) \\
\mu = (8g+4)\lambda - g \delta_{irr} - 
\sum_{i=1}^{g-1} 4i(g-i)\delta_{i} \\
\theta_1 =  4g(g-1)\psi_1 - 12 \lambda + 
\delta_{irr} - \sum_{i=1}^{g-1}
4i(i-1) \delta_{i}.
\end{cases}
\]
Then, we
have the following:
\begin{enumerate}
\renewcommand{\labelenumi}{(\arabic{enumi})}
\item
$\mu$ and $\theta_1$ are weakly positive over any finite subsets of $\MM_{g,1}$.
In particular, 
\[
\QQ_{+} \mu + \QQ_{+}\theta_1 + \QQ_{+} \delta_{irr} +
\sum_{i=1}^{g-1} \QQ_{+} \delta_{i}
\subseteq \Nef(\bMM_{g,1}; \MM_{g,1}),
\]
where $\QQ_{+} = \{ x \in \QQ \mid x \geq 0 \}$.
\rom{(}Note that $\mu = \theta_1 = 0$ if $g=1$, and
$\mu=0$ if $g=2$.\rom{)}

\item
We assume $g = 1$.
Then, $a \mu + b \theta_1 + c_{irr} \delta_{irr}$ is nef
over $\MM_{1,1}$ if and only if $c_{irr} \geq 0$.

\item
We assume $g \geq 2$.
If a $\QQ$-divisor 
\[
D = a \mu + b \theta_1 + c_{irr} \delta_{irr}
+ \sum_{i=1}^{g-1} c_{i}\delta_{i}
\]
is nef over $\MM_{g,1}$, then
$b, c_{irr}, c_{1}, \ldots, c_{g-1}$ are
non-negative.
\end{enumerate}
\end{Corollary}

\Proof
(1) $\mu$ is weakly positive over any finite subsets of $\MM_{g,1}$
by \cite[Theorem~B]{MoRB} or Remark~\ref{rem:slope:inq:mg}, and
(2) of Proposition~\ref{prop:samp:pamp:pullback}.
Moreover, $\theta_1$ is weakly positive over any finite subsets of $\MM_{g,1}$
by virtue of the case $T = L= \{1\}$ in 
Theorem~\ref{thm:slope:inq:m:g:n}.

\medskip
(2)
This is obvious because $\mu = \theta_1 = 0$.

\medskip
(3) We assume that $D$ is nef over $\MM_{g,1}$.
Let $C$ be a smooth curve of genus $g$, and 
$\Delta$ the diagonal of $C \times C$.
Let $p : C \times C \to C$ be the projection to the first factor.
Then, $\Delta$ gives rise to a section of $p$.
Hence, we get a morphism $\varphi_1 : C \to \bMM_{g,1}$ with
$\varphi_1(C) \subseteq \MM_{g,1}$. By our assumption, 
$\deg(\varphi_1^*(D)) \geq 0$.
On the other hand,
\[
\deg(\varphi_1^*(\mu)) = \deg(\varphi_1^*(\delta_{irr})) = 
\deg(\varphi_1^*(\delta_{1})) =
\cdots = \deg(\varphi_1^*(\delta_{g-1})) = 0
\]
and $\deg(\varphi_1^*(\theta_1))  = 8g(g-1)^2$. 
Thus, $b \geq 0$.

Let $f_2 : X_2 \to Y_2$ be a hyperelliptic fibered surface and
$\Gamma_2$ a section as in  Proposition~\ref{prop:irr:one} for $i=0$.
Let $\varphi_2 : Y_2 \to \bMM_{g,1}$ be the induced morphism.
Then, $\varphi_2(Y_2) \cap \MM_{g,1} \not= \emptyset$,
\[
\deg(\varphi_2^*(\mu)) = \deg(\varphi_2^*(\delta_{1})) = \cdots = 
\deg(\varphi_2^*(\delta_{g-1})) = 0
\]
and $\deg(\varphi_2^*(\delta_{irr})) = \deg(\delta_{irr}(X_2/Y_2)) > 0$.
On the generic fiber, $\Gamma_2$ is a ramification point of
the hyperelliptic covering. Thus,
\[
L_2 = \omega_{X_2/Y_2} \otimes \OO_{X_2}(-(2g-2)\Gamma_2)
\]
satisfies the conditions of Lemma~\ref{lem:self:2:zero}.
Thus, $(L_2^2) = 0$, which says us that $\deg(\varphi^*_2(\theta_1)) = 0$.
Therefore, we get $c_{irr} \geq 0$.

Finally we fix $i$ with $1 \leq i \leq g-1$.
Let $f_3 : X_3 \to Y_3$ be a hyperelliptic fibered surface and $\Gamma_3$
a section as in Proposition~\ref{prop:a:ga:one}.
Let $\varphi_3 : Y_3 \to \bMM_{g,1}$ be the induced morphism.
Then, $\varphi_3(Y_3) \cap \MM_{g,1} \not= \emptyset$,
$\deg(\varphi_3^*(\mu)) = 0$,
$\deg(\varphi_3^*(\delta_{l})) = 0$ ($l \not= i$)
and $\deg(\varphi_3^*(\delta_{i})) = \deg(\delta_{i}(X_3/Y_3)) > 0$.
Let $\Sigma_3$ be the set of critical values of $f_3$.
For each $P \in \Sigma_3$, let $E_P$ 
be the component of genus $i$ in $f_3^{-1}(P)$.
On the generic fiber, $\Gamma_2$ is a ramification point of
the hyperelliptic covering. Thus,
\[
L_3 = \omega_{X_3/Y_3} \otimes \OO_{X_3}\left(-(2g-2)\Gamma_3 + 
\sum_{P \in \Sigma_3} (2i-1) E_P\right)
\]
satisfies the conditions of Lemma~\ref{lem:self:2:zero}.
Therefore, $(L_3^2) = 0$, which says us that 
$\deg(\varphi^*_3(\theta_1)) = 0$. Hence, we get 
$c_{i} \geq 0$.
\QED

\begin{Corollary}[$\ch(k) = 0$]
\label{cor:slope:cone:mg2}
Let $\bMM_{g,2} = \bMM_{g,\{1,2\}}$ be the moduli space of two-pointed
stable curves of genus $g \geq 2$.
We set $\delta_{i}, \sigma_i,
\mu, \theta_{1,2} \in \Pic(\bMM_{g,2}) \otimes \QQ$ as follows:
\[
\begin{cases}
\delta_i = \delta_{\{(i, \emptyset), (g-i, \{1, 2\})\}}
\quad(1 \leq i \leq g) \\
\sigma_i = \delta_{\{(i, \{1\}), (g-i, \{2\})\}}
\quad(1 \leq i \leq g-1) \\
\mu = (8g+4)\lambda - g\delta_{irr} -
\sum_{i=1}^{g-1} 4i(g-i)
\sigma_{i} - \sum_{i=1}^g 4i(g-i) \delta_{i}, \\[6pt]
\theta_{1,2} = 
(g-1)(g+1)(\psi_1 + \psi_2) - 12 \lambda + \delta_{irr} \\[0pt]
\qquad\qquad -\sum_{i=1}^{g-1} (2i-g-1)(2i-g+1) \sigma_{i}
-\sum_{i=1}^{g}4i(i-1) \delta_{i}.
\end{cases}
\]
Then, we have the following.
\begin{enumerate}
\renewcommand{\labelenumi}{(\arabic{enumi})}
\item
$\mu$ and $\theta_{1,2}$ are weakly positive over
any finite subsets of $\MM_{g,2}$.
In particular,
\[
\QQ_{+} \mu + \QQ_{+} \theta_{1,2} + \QQ_{+} \delta_{irr} +
\sum_{i=1}^{g-1}
\QQ_{+} \sigma_{i} +
\sum_{i=1}^g \QQ_{+} \delta_{i}
\subseteq \Nef(\bMM_{g,2}; \MM_{g,2}).
\]

\item
If a $\QQ$-divisor
\[
D = a \mu + b \theta_{1,2} + c_{irr} \delta_{irr} +
\sum_{i=1}^{g-1}
c_{i}\sigma_{i} +
\sum_{i=1}^g d_{i}\delta_{i}
\]
on $\bMM_{g,2}$ is
nef over $\MM_{g,2}$, then
\[
b \geq 0,\  c_{irr} \geq 0,\ 
c_{i} \geq 0\ (\forall i=1,\ldots,g-1),\  
d_{i} \geq 0\ (\forall i=1, \ldots, g).
\]

\item
Here we set $\sigma$, $\mu'$ and $\theta'_{1,2}$
as follows:
\[
\begin{cases}
\sigma = \delta_{irr} + 
\sum_{i=1}^{g-1}
\sigma_{i}, \\
\mu' = (8g+4) \lambda - g \sigma - \sum_{i=1}^g 4i(g-i)
\delta_{i}, \\
\theta'_{1,2} =
(g-1)(g+1)(\psi_1 + \psi_2) - 12 \lambda + \sigma 
-\sum_{i=1}^{g}4i(i-1) \delta_{i}.
\end{cases}
\]
Then, we have
\[
\QQ_{+} \mu' + \QQ_{+} \theta'_{1,2} + \QQ_{+} \sigma +
\sum_{i=1}^g \QQ_{+} \delta_{i}
\subseteq \Nef(\bMM_{g,2}; \MM_{g,2}).
\]
Moreover,
if a $\QQ$-divisor $a \mu' + b \theta'_{1,2} + c
\sigma + \sum_{i=1}^{g} d_{i} \delta_{i}$ on 
$\bMM_{g,2}$ is nef over $\MM_{g,2}$,
then
$b$, $c$, $d_{1}, \ldots, d_{g}$ are non-negative.
\end{enumerate}
\end{Corollary}

\Proof
(1) By \cite[Theorem~B]{MoRB} or Remark~\ref{rem:slope:inq:mg},
and (2) of Proposition~\ref{prop:samp:pamp:pullback},
$\mu$ is weakly positive over any finite subsets of $\MM_{g,2}$.
Further, $\theta_{1,2}$ is weakly positive over 
any finite subsets of $\MM_{g,2}$
by the case $T = L= \{1, 2\}$ in 
Theorem~\ref{thm:slope:inq:m:g:n}.

\medskip
(2) We assume that $D$ is nef over $\MM_{g,2}$.
Let $C$ be a smooth curve of genus $g$, and $\Delta$ 
the diagonal of $C \times C$.
Let $p : C \times C \to C$ and $q : C \times C \to C$ be 
the projection to the first factor
and the second factor respectively. Moreover, let $\vartheta$ be a line bundle
on $C$ with $\vartheta^{\otimes 2} = \omega_C$ and
$L_n = p^*(\vartheta^{\otimes n}) \otimes q^*(\vartheta^{\otimes n}) 
\otimes \OO_{C \times C}((n-1)\Delta)$.
For $n \geq 3$, let $T_n$ be a general member of
$\left| L_n \right|$.
Then, since $(L_n^2) > 0$,
by Lemma~\ref{lem:base:pt:free},
$T_n$ is smooth and irreducible.
Moreover, $T_n$ meets $\Delta$ transversally.
Then, we have two morphisms $p_n : T_n \to C$ and $q_n : T_n \to C$
given by $T_n \hookrightarrow C \times C \overset{p}{\longrightarrow} C$
and $T_n \hookrightarrow C \times C \overset{q}{\longrightarrow} C$
respectively.
Let $\Gamma_{p_n}$ and $\Gamma_{q_n}$ be the graph of $p_n$ and $q_n$
in $C \times T_n$ respectively.
Then, it is easy to see that
$\Gamma_{p_n}$ and $\Gamma_{q_n}$ meet transversally, and
$(\Gamma_{p_n} \cdot \Gamma_{q_n}) = (T_n \cdot \Delta) = 2g-2$.
Let $X \to C \times T_n$ be the blowing-ups at points in 
$\Gamma_{p_n} \cap \Gamma_{q_n}$, and let
$\overline{\Gamma}_{p_n}$ and $\overline{\Gamma}_{q_n}$
be the strict transform of $\Gamma_{p_n}$ and $\Gamma_{q_n}$
respectively.
Then, $\overline{\Gamma}_{p_n}$ and $\overline{\Gamma}_{q_n}$
give rise to two non-crossing sections of $X \to T_n$.
Moreover,
\[
(\omega_{X/T_n} \cdot \overline{\Gamma}_{p_n}) = 
(\omega_{C \times T_n/T_n} \cdot \Gamma_{p_n}) = 
2(g-1)\deg(\Gamma_{p_n} \to C)
= 2(g-1)(ng - 1).
\]
In the same way, $(\omega_{X/T_n} \cdot \overline{\Gamma}_{q_n}) = 
2(g-1)(ng-1)$.
Let $\pi_n : T_n \to \bMM_{g,2}$ be the induced morphism.
Then, we can see that
$\deg(\pi_n^*(\lambda)) =  \deg(\pi_n^*(\sigma_{i}))
=\deg(\pi_n^*(\delta_{i})) = 0$
for all $i = 1,\ldots,g-1$.
Moreover, $\deg(\pi_n^*(\psi_1 + \psi_2)) = 4(g-1)(ng-1)$, and
$\deg(\pi_n^*(\delta_{g})) = 2(g-1)$.
Thus,
\[
\deg(\pi_n^*(D)) = 4(g+1)(g-1)^2(ng-1)b - 8g(g-1)^2d_{g} \geq 0
\]
for all $n \geq 3$. Therefore, we get $b \geq 0$.

Let $f_2 : X_2 \to Y_2$ be a hyperelliptic fibered surface and $\Gamma_2$
a section as in Proposition~\ref{prop:irr:two} for $i=0$.
Then, $\Gamma_2$ and $j(\Gamma_2)$ gives two points of $X_2$ over $Y_2$.
Let $\varphi_2 : Y_2 \to \bMM_{g,2}$ be the induced morphism.
Then, $\varphi_2(Y_2) \cap \MM_{g,2} \not= \emptyset$,
$\deg(\varphi_2^*(\mu)) = 0$,
$\deg(\varphi_2^*(\sigma_{i})) = 0$ for all
$i = 1,\ldots,g-1$, and
$\deg(\varphi_2^*(\delta_i)) = 0$ for all $i=1,\ldots,g$.
Moreover, $\deg(\varphi_2(\delta_{irr}))  > 0$.
On the generic fiber, two points arising from $\Gamma_2$ and $j(\Gamma_2)$
are invariant under the action of the hyperelliptic involution.
Thus,
\[
L_2 = \omega_{X_2/Y_2} \otimes \OO_{X_2}(-(g-1)(\Gamma_2+j(\Gamma_2))
\]
satisfies the conditions of Lemma~\ref{lem:self:2:zero}.
Thus, $(L_2^2) = 0$, which says us that 
$\deg(\varphi^*_2(\theta_{1,2})) = 0$. 
Thus, we get $c_{irr} \geq 0$.

We fix $i$ with $1 \leq i \leq g$.
Let $f_3 : X_3 \to Y_3$ be a hyperelliptic fibered surface and $\Gamma_3$
a section as in Proposition~\ref{prop:a:ga:two}.
Let $\varphi_3 : Y_3 \to \bMM_{g,2}$ be the induced morphism
arising from the $2$-pointed curve 
$\{ f_3: X_3 \to Y_3; \Gamma_3, j(\Gamma_3) \}$.
Then, $\varphi_3(Y_3) \cap \MM_{g,2} \not= \emptyset$,
$\deg(\varphi_3^*(\mu)) = 0$, 
$\deg(\varphi_3^*(\sigma_{s})) = 0$ ($\forall
s=1,\ldots,g-1$),
$\deg(\varphi_3^*(\delta_{s})) = 0$ ($\forall s \not= i$)
and $\deg(\varphi_3^*(\delta_{i})) = \deg(\delta_i(X_3/Y_3)) > 0$.
Let $\Sigma_3$ be the set of critical values of $f_3$.
For each $P \in \Sigma_3$, let $E_P$ be 
the component of genus $i$ in $f_3^{-1}(P)$.
On the generic fiber, two points arising from $\Gamma_2$ and $j(\Gamma_2)$
are invariant under the action of the hyperelliptic involution.
Thus,
\[
L_3 = \omega_{X_3/Y_3} \otimes \OO_{X_3}\left(-(g-1)(\Gamma_3 + j(\Gamma_3)) 
+ \sum_{P \in \Sigma_3} (2i-1) E_P\right)
\]
satisfies the conditions of Lemma~\ref{lem:self:2:zero}.
Therefore, $(L_3^2) = 0$, which says us that
$\deg(\varphi^*_3(\theta_{1,2})) = 0$. Hence, we get
$d_{i} \geq 0$.

Finally we fix $i$ with $1 \leq i \leq g-1$.
Let $f_4 : X_4 \to Y_4$ be a hyperelliptic fibered surface and 
$\Gamma_4, \Gamma'_4$
sections as in Proposition~\ref{prop:g:a:two:section}.
Let $\varphi_4 : Y_4 \to \bMM_{g,2}$ be the induced morphism.
Then, $\varphi_4(Y_4) \cap \MM_{g,2} \not= \emptyset$,
$\deg(\varphi_4^*(\mu)) = 0$,
$\deg(\varphi_4^*(\delta_{s})) = 0$ ($\forall s$),
$\deg(\varphi_4^*(\sigma_{s})) = 0$ ($\forall
s \not= i$), and
$\deg(\varphi_4^*(\sigma_{i})) > 0$.
Let $\Sigma_4$ be the set of critical values of $f_4$.
For each $P \in \Sigma_4$, let $E_P$ be the component of genus $i$ in
$f_4^{-1}(P)$. On the generic fiber, $\Gamma_4$ and $\Gamma'_4$ are a
ramification point of the hyperelliptic covering. Thus,
\[
L_4 = \omega_{X_4/Y_4} \otimes \OO_{X_4}\left(-(g-1)(\Gamma_4 + \Gamma'_4) + 
\sum_{P \in \Sigma_4} \left( (2i-1) - (g-1) \right) E_P\right)
\]
satisfies the conditions of Lemma~\ref{lem:self:2:zero}.
Therefore, $(L_4^2) = 0$, which says us that
$\deg(\varphi^*_4(\theta_{1,2})) = 0$. Hence, we get 
$c_{i} \geq 0$.

\medskip
(3) There are non-negative integers
$e_{i}$ and $f_{i}$
($1 \leq i \leq g-1$) with
\[
\mu' = \mu + \sum_{i=1}^{g-1} e_{i}\sigma_{i}
\quad\text{and}\quad
\theta'_{1,2} = \theta_{1,2} + \sum_{i=1}^{g-1} f_{i}\sigma_{i}.
\]
Thus, (3) is a consequence of (1) and (2).
\QED

\section{The proof of the main result}
Throughout this section,
we fix an integer $g \geq 3$.
The purpose of this section is to prove the following theorem.

\begin{Theorem}[$\ch(k) = 0$]
\label{thm:cone:mg:codim:two}
A $\QQ$-divisor $a \mu + b_{irr} \delta_{irr} +
\sum_{i=1}^{[g/2]} b_{i}
\delta_{i}$ on $\bMM_{g}$
is nef over $\bMM_{g}^{[1]}$ if and only if
the following system of inequalities hold:
\begin{align*}
& a \geq \max \left\{ \frac{b_{i}}{4i(g-i)} \mid 
i = 1, \ldots, [g/2] \right\} \\
& B_0 \geq B_1 \geq B_2 \geq \cdots \geq B_{[g/2]}, \\
& B^*_{[g/2]} \geq \cdots \geq B^*_2 \geq B^*_1 \geq B^*_0,
\end{align*}
where $B_0$, $B_0^*$, $B_i$ and $B_i^*$ 
\rom{(}$i=1, \ldots, [g/2]$\rom{)} are given by
\[
B_0 = 4b_{irr},\quad
B^*_0 = \frac{4b_{irr}}{g(2g-1)},\quad
B_i = \frac{b_i}{i(2i+1)}\quad\text{and}\quad
B^*_i = \frac{b_i}{(g-i)(2(g-i)+1)}.
\]
\end{Theorem}

\Proof
In the following proof, we denote $\delta_{i}$ by $\delta_{\{i,g-i\}}$.
Moreover, we set
\[
\overline{\upsilon}_g =
\{ \{i,j\} \mid 1 \leq i, j \leq g,\ i+j=g \}.
\]

For a $\QQ$-divisor $D = a \mu + b_{irr} \delta_{irr} +
\sum_{\{i,j\} \in \overline{\upsilon}_g} b_{\{i,j\}}
\delta_{\{i,j\}}$, let us consider the following inequalities:
{\allowdisplaybreaks
\begin{align}
\addtocounter{Claim}{1}
\label{eqn:thm:cone:mg:codim:two:A}
& a \geq \frac{b_{\{s,t\}}}{4st}\quad 
(\forall \{s,t\} \in \overline{\upsilon}_g) \\
\addtocounter{Claim}{1}
\label{eqn:thm:cone:mg:codim:two:B}
& 4b_{irr} \geq \frac{b_{\{s,t\}}}{s(2s+1)},\quad
  \frac{b_{\{s,t\}}}{t(2t+1)} \geq \frac{4b_{irr}}{g(2g-1)}
\quad\text{\rom{(}$\forall \{s,t\} \in \overline{\upsilon}_g$ with 
$s\leq t$\rom{)}}
\\
\addtocounter{Claim}{1}
\label{eqn:thm:cone:mg:codim:two:C}
& \frac{b_{\{l,k\}}}{l(2l+1)} \geq
\frac{b_{\{s,t\}}}{s(2s+1)},\quad 
\frac{b_{\{s,t\}}}{t(2t+1)} \geq
\frac{b_{\{l,k\}}}{k(2k+1)} \\[-8pt]
& \phantom{\frac{b_{\{l,k\}}}{l(2l+1)} \geq
\frac{b_{\{s,t\}}}{s(2s+1)} \geq \frac{b_{\{s,t\}}}{t(2t+1)}}
\text{\rom{(}$\forall \{s, t\}, \{l,k\} \in
\overline{\upsilon}_g$ with $l < s \leq t < k$\rom{)}} \notag \\
\addtocounter{Claim}{1}
\label{eqn:thm:cone:mg:codim:two:D}
& a \geq 0,\quad b_{irr} \geq 0,\quad
b_{\{s,t\}} \geq 0\quad (\forall \{s,t\} \in \overline{\upsilon}_g)
\end{align}}
Let
$\beta : \bMM_{g-1,2} \to \bMM_g$ and
$\alpha_{s,t} : \bMM_{s,1} \times \bMM_{t,1} \to \bMM_g$
($\{s,t\} \in \overline{\upsilon}_g$)
be the clutching maps.
First, we claim the following.

\begin{Claim}
\label{thm:cone:mg:codim:two:claim:1}
The following are equivalent.
\begin{enumerate}
\renewcommand{\labelenumi}{(\arabic{enumi})}
\item 
$\beta^*(D)$ is nef over $\MM_{g-1,2}$ and
$\alpha_{s,t}^*(D)$ is nef over $\MM_{s,1} \times \MM_{t,1}$
for all $\{ s, t \} \in \overline{\upsilon}_g$

\item
\eqref{eqn:thm:cone:mg:codim:two:A},
\eqref{eqn:thm:cone:mg:codim:two:B},
\eqref{eqn:thm:cone:mg:codim:two:C} and
\eqref{eqn:thm:cone:mg:codim:two:D} hold.
\end{enumerate}
\end{Claim}

On $\bMM_{g-1,2}$, we define $\sigma$ and $\delta_i$ ($i=1,\ldots,g-1$)
as in Corollary~\ref{cor:slope:cone:mg2}.
Moreover, we set
\begin{align*}
\mu' & = (8g-4)\lambda -(g-1)\sigma -
\sum_{i=1}^{g-1}4i(g-1-i)\delta_{i}, \\
\theta' & = (g-2)g(\psi_1 + \psi_2) - 12\lambda + \sigma
- \sum_{i=1}^{g-1}4i(i-1)\delta_{i}.
\end{align*}
Then, by using (2) of Corollary~\ref{cor:formula:clutching:map},
we can see 
\begin{multline*}
\beta^*(D) =
\frac{(g-1)(g-2)(2g-1)a - 3b_{irr}}{g(g-2)(2g-1)}\mu' +
\frac{ag-b_{irr}}{g(g-2)} \theta' \\ +
\frac{(g-1)(2g+1)b_{irr}}{g(2g-1)}\sigma +
\sum_{i=1}^{g-1}\left(b_{\{i,g-i\}} - \frac{4i(2i+1)}{g(2g-1)}b_{irr}
\right)\delta_{i}.
\end{multline*}
Thus,
by Corollary~\ref{cor:slope:cone:mg2},
if $\beta^*(D)$ is nef over $\bMM_{g-1,2}$, then
{\allowdisplaybreaks
\begin{align}
\addtocounter{Claim}{1}
\label{eqn:thm:cone:mg:codim:two:1}
& ag \geq b_{irr} \geq 0 \\
\addtocounter{Claim}{1}
\label{eqn:thm:cone:mg:codim:two:2}
& b_{\{i,g-i\}} \geq \frac{4i(2i+1)}{g(2g-1)}b_{irr},
\quad(i = 1, \ldots, g-1)
\end{align}}

\medskip
Here we set $\mu'_1 = \theta'_1 = 0$ on $\bMM_{1,1}$, and
\begin{align*}
\mu'_e & = \frac{1}{e-1}
\left((8e+4)\lambda - e\delta_{irr} -
\sum_{l=1}^{e-1} 4l(e-l) \delta_{l}\right), \\
\theta'_e & = \frac{1}{e-1}\left(
4e(e-1)\psi_1 - 12 \lambda + \delta_{irr} - \sum_{l=1}^{e-1}
4l(l-1)\delta_{l}\right)
\end{align*}
on $\bMM_{e,1}$ ($e \geq 2$), where
$\delta_l$'s are defined as in Corollary~\ref{cor:slope:cone:mg1}.
Let us fix $\{s,t\} \in \overline{\upsilon}_g$.
Then, by using (1) of Corollary~\ref{cor:formula:clutching:map},
we can see
\[
\alpha_{s,t}^*(D) = p^*(D_s) + q^*(D_t),
\]
where $p : \bMM_{s,1} \times \bMM_{t,1}
\to \bMM_{s,1}$ and $q: \bMM_{s,1} \times \bMM_{t,1}
\to \bMM_{t,1}$ are the projections, and
$D_s \in \Pic(\bMM_{s,1}) \otimes \QQ$ and
$D_t \in \Pic(\bMM_{t,1}) \otimes \QQ$ are given by
\begin{multline*}
D_s = \frac{4(g-1)s(2s+1)a - 3b_{\{s,t\}}}{4s(2s+1)}\mu'_s
+ \frac{4sta - b_{\{s,t\}}}{4s}\theta'_s \\
+ \left(b_{irr} - \frac{b_{\{s,t\}}}{4s(2s+1)} \right)\delta_{irr} +
\sum_{l=1}^{s-1}\left( b_{\{l,g-l\}} -
\frac{l(2l+1)}{s(2s+1)}b_{\{s,t\}}\right) \delta_{l}
\end{multline*}
and
\begin{multline*}
D_t = \frac{4(g-1)t(2t+1)a - 3b_{\{s,t\}}}{4t(2t+1)}\mu'_t
+ \frac{4sta - b_{\{s,t\}}}{4t}\theta'_s \\
+ \left(b_{irr} - \frac{b_{\{s,t\}}}{4t(2t+1)} \right)\delta_{irr} +
\sum_{l=1}^{t-1}\left( b_{\{l,g-l\}} -
\frac{l(2l+1)}{t(2t+1)}b_{\{s,t\}}\right) \delta_{l}.
\end{multline*}
Thus, by using Corollary~\ref{cor:slope:cone:mg1} and
Lemma~\ref{lem:nef:over:product}, if
$\alpha_{s,t}^*(D)$ is nef over
$\MM_{s,1} \times \MM_{t,1}$, then
{
\allowdisplaybreaks
\begin{align}
\addtocounter{Claim}{1}
\label{eqn:thm:cone:mg:codim:two:3}
& 4sta \geq b_{\{s,t\}}  \\
\addtocounter{Claim}{1}
\label{eqn:thm:cone:mg:codim:two:4}
& b_{irr} \geq \frac{b_{\{s,t\}}}{4s(2s+1)},\quad
  b_{irr} \geq \frac{b_{\{s,t\}}}{4t(2t+1)} \\
\addtocounter{Claim}{1}
\label{eqn:thm:cone:mg:codim:two:5}
& b_{\{l,g-l\}} \geq
\frac{l(2l+1)}{s(2s+1)}b_{\{s,t\}}
\quad (l=1, \ldots, s-1) \\
\addtocounter{Claim}{1}
\label{eqn:thm:cone:mg:codim:two:6}
& b_{\{l,g-l\}} \geq
\frac{l(2l+1)}{t(2t+1)}b_{\{s,t\}}
\quad (l=1, \ldots, t-1)
\end{align}
}

Therefore, (1) implies
\eqref{eqn:thm:cone:mg:codim:two:1} -- \eqref{eqn:thm:cone:mg:codim:two:6}.
Conversely, we assume
\eqref{eqn:thm:cone:mg:codim:two:1} -- \eqref{eqn:thm:cone:mg:codim:two:6}.
Then by using\eqref{eqn:thm:cone:mg:codim:two:1} and \eqref{eqn:thm:cone:mg:codim:two:2},
we can see \eqref{eqn:thm:cone:mg:codim:two:D}.
Thus, we have
\begin{align*}
ag-b_{irr} \geq 0 & \quad\Longrightarrow\quad
(g-1)(g-2)(2g-1)a - 3b_{irr} \geq 0 \\
4sta \geq b_{\{s,t\}} &\quad\Longrightarrow\quad
\text{$4(g-1)s(2s+1)a \geq  3b_{\{s,t\}}$ and
$4(g-1)t(2t+1)a \geq 3b_{\{s,t\}}$}.
\end{align*}
Therefore, by Corollary~\ref{cor:slope:cone:mg1}, Corollary~\ref{cor:slope:cone:mg2}
and Lemma~\ref{lem:nef:over:product},
we can see that 
$\beta^*(D)$ is nef over $\MM_{g-1,2}$ and
$\alpha_{s,t}^*(D)$ is nef over $\MM_{s,1} \times \MM_{t,1}$
for all $\{ s, t \} \in \overline{\upsilon}_g$.
Hence it is sufficient to see that
the system of inequalities
\eqref{eqn:thm:cone:mg:codim:two:1} -- \eqref{eqn:thm:cone:mg:codim:two:6}
is equivalent to
\eqref{eqn:thm:cone:mg:codim:two:A} -- \eqref{eqn:thm:cone:mg:codim:two:C}
under the assumption \eqref{eqn:thm:cone:mg:codim:two:D}.

\medskip
The case $s=1,t=g-1$ in \eqref{eqn:thm:cone:mg:codim:two:3}
and the case $i=g-1$ in \eqref{eqn:thm:cone:mg:codim:two:2}
produce inequalities
\[
4(g-1)a \geq b_{\{1,g-1\}}\quad\text{and}\quad
b_{\{1,g-1\}} \geq \frac{4(g-1)}{g} b_{irr}
\]
respectively,
which gives rise to
\eqref{eqn:thm:cone:mg:codim:two:1}.
Moreover, it is easy to see that
\eqref{eqn:thm:cone:mg:codim:two:2} and
\eqref{eqn:thm:cone:mg:codim:two:4} are
equivalent to
\eqref{eqn:thm:cone:mg:codim:two:B}, so that
it is sufficient to see that
\eqref{eqn:thm:cone:mg:codim:two:5} and
\eqref{eqn:thm:cone:mg:codim:two:6} are equivalent to
\eqref{eqn:thm:cone:mg:codim:two:C}.

From now on , we assume $s \leq t$.
Since $s(2s+1) \leq t(2t+1)$,
\eqref{eqn:thm:cone:mg:codim:two:5} and
\eqref{eqn:thm:cone:mg:codim:two:6} are equivalent to saying that
{
\allowdisplaybreaks
\begin{align}
\addtocounter{Claim}{1}
\label{eqn:thm:cone:mg:codim:two:7}
& \frac{b_{\{l,k\}}}{l(2l+1)} \geq
\frac{b_{\{s,t\}}}{s(2s+1)}
\quad(1 \leq l < s) \\
\addtocounter{Claim}{1}
\label{eqn:thm:cone:mg:codim:two:8}
&\frac{b_{\{l,k\}}}{l(2l+1)} \geq
\frac{b_{\{s,t\}}}{t(2t+1)}
\quad(s < l < t),
\end{align}
}
where $k = g-l$. 
In \eqref{eqn:thm:cone:mg:codim:two:7},
$t< k \leq g-1$, Thus, \eqref{eqn:thm:cone:mg:codim:two:7}
is nothing more than
\[
\frac{b_{\{l,k\}}}{l(2l+1)} \geq
\frac{b_{\{s,t\}}}{s(2s+1)}
\quad(1 \leq l < s \leq t < k \leq g-1)
\]
Moreover, in \eqref{eqn:thm:cone:mg:codim:two:8},
$s < k < t$.
Thus, \eqref{eqn:thm:cone:mg:codim:two:8} is nothing more that
\[
\frac{b_{\{l,k\}}}{k(2k+1)} \geq
\frac{b_{\{s,t\}}}{t(2t+1)}\quad
(1 \leq s < l \leq k < t \leq g-1).
\]
Therefore, replacing $\{s,t\}$ and $\{l,k\}$, we have
\[
\frac{b_{\{s,t\}}}{t(2t+1)} \geq
\frac{b_{\{l,k\}}}{k(2k+1)}\quad
(1 \leq l < s \leq t < k \leq g-1).
\]
Thus, we get Claim~\ref{thm:cone:mg:codim:two:claim:1}.

\medskip
By Claim~\ref{thm:cone:mg:codim:two:claim:1},
it is sufficient to show the following claim
to complete the proof of Theorem~\ref{thm:cone:mg:codim:two}.

\begin{Claim}
\begin{enumerate}
\renewcommand{\labelenumi}{(\arabic{enumi})}
\item
$D$ is nef over $\bMM_{g}^{[1]}$ if and only if
$D$ is nef over $\MM_g$, $\beta^*(D)$ is nef over $\MM_{g-1,2}$, and
$\alpha_{s,t}^*(D)$ is nef over $\MM_{s,1} \times \MM_{t,1}$
for all $\{ s, t \} \in \overline{\upsilon}_g$.

\item
$D$ is nef over $\MM_{g}$ if and only if
\eqref{eqn:thm:cone:mg:codim:two:D} holds.

\item
\eqref{eqn:thm:cone:mg:codim:two:A},
\eqref{eqn:thm:cone:mg:codim:two:B} and
\eqref{eqn:thm:cone:mg:codim:two:C} imply
\eqref{eqn:thm:cone:mg:codim:two:D}
\end{enumerate}
\end{Claim}

(1) is obvious because
\[
\bMM_{g}^{[1]} = \MM_g \cup \beta_g(\MM_{g-1,2}) \cup 
\bigcup_{\{s,t\} \in \overline{\upsilon}_g}
\alpha_{s,t}(\MM_{s,1} \times \MM_{t,1}).
\]
(2) is a consequence of \cite[Theorem~C]{MoRB}.
For (3), let us consider the case $s=1$, $t=g-1$ in
\eqref{eqn:thm:cone:mg:codim:two:B}. Then, we have
\[
12b_{irr} \geq b_{\{1,g-1\}}\quad\text{and}\quad
b_{\{1,g-1\}} \geq \frac{4(g-1)}{g}b_{irr},
\]
which imply $b_{irr} \geq 0$. Thus, we can see
\eqref{eqn:thm:cone:mg:codim:two:D}
using \eqref{eqn:thm:cone:mg:codim:two:A} and
\eqref{eqn:thm:cone:mg:codim:two:B}.
\QED

\begin{Corollary}[$\ch(k) = 0$]
\label{cor:cone:mg:codim:two}
Let $\widetilde{\Delta}_{irr}$ and
$\widetilde{\Delta}_i$ \rom{(}$i=1,\ldots,[g/2]$\rom{)}
be the normalizations of
the boundary components $\Delta_{irr}$ and $\Delta_{i}$ on $\bMM_g$, and
$\rho_{irr} : \widetilde{\Delta}_{irr} \to \bMM_g$
and $\rho_i : \widetilde{\Delta}_{i} \to \bMM_g$
the induced morphisms.
Then, a $\QQ$-divisor $D$ on $\bMM_g$ is nef over $\bMM_g^{[1]}$ if and only if
the following are satisfied:
\begin{enumerate}
\renewcommand{\labelenumi}{(\arabic{enumi})}
\item
$D$ is weakly positive at any points of $\MM_g$.

\item
$\rho_{irr}^*(D)$ is weakly positive at any points of 
$\rho_{irr}^{-1}(\bMM_g^{[1]})$.

\item
$\rho_{i}^*(D)$ is weakly positive at any points of 
$\rho_{i}^{-1}(\bMM_g^{[1]})$ for all $i$.
\end{enumerate}
\end{Corollary}

\Proof
Let $\beta : \bMM_{g-1,2} \to \bMM_g$ be the clutching map.
Then, there is a finite and
surjective morphism $\beta' : \bMM_{g-1,2} \to \widetilde{\Delta}_{irr}$
with $\beta = \rho_{irr} \cdot \beta'$.
Further, for $1 \leq i \leq [g/2]$,
let $\alpha_{i,g-i} : \bMM_{i,1} \times \bMM_{g-i,1} \to \bMM_g$
be the clutching map. Then, there is 
a finite and surjective morphism $\alpha'_{i,g-i}:  \bMM_{i,1} \times \bMM_{g-i,1} \to 
\widetilde{\Delta}_{i}$
with $\alpha_{i,g-i} = \rho_{i} \cdot \alpha'_{i,g-i}$.
In particular, $\widetilde{\Delta}_{irr}$ and
$\widetilde{\Delta}_{i}$'s are $\QQ$-factorial.
Therefore, if $D$ satisfies (1), (2) and (3), then $D$ is nef over
$\bMM_g^{[1]}$. 

Conversely, we assume that $D$ is nef over $\bMM_g^{[1]}$.
(1) is nothing more than \cite[Theorem~C]{MoRB}.
As in Theorem~\ref{thm:cone:mg:codim:two}, we set
$D = a \mu + b_{irr} \delta_{irr} +
\sum_{i=1}^{[g/2]} b_{i}
\delta_{i}$ on $\bMM_{g}$.
If we trace-back the proof of Theorem~\ref{thm:cone:mg:codim:two},
we can see that
\[
\beta^*(D) \in \QQ_{+}\mu' +
\QQ_{+} \theta' + \QQ_{+} \sigma + \sum_i \QQ_{+} \delta_{i}.
\]
Here $\mu'$ and $\theta'$ are weakly positive at any points
of $\MM_{g-1,2}$ by (1) of Corollary~\ref{cor:slope:cone:mg2}.
Thus, so is $\beta^*(D) = {\beta'}^*(\rho_{irr}^*(D))$.
Therefore, by virtue of (2) of Proposition~\ref{prop:samp:pamp:push},
$\beta'_*(\beta^*(D)) = \deg(\beta')\rho_{irr}^*(D)$ is
weakly positive at any points of 
$\rho_{irr}^{-1}(\bMM_g^{[1]})$.
Finally, let us consider (3).
As in the proof of Theorem~\ref{thm:cone:mg:codim:two},
there are $D_i \in \Pic(\bMM_{i,1}) \otimes \QQ$ and
$D_{g-i} \in  \Pic(\bMM_{g-i,1}) \otimes \QQ$ with
$\alpha_{i,g-i}^*(D) = p^*(D_i) + q^*(D_{g-i})$, where
$p : \bMM_{i,1} \times \bMM_{g-i,1} \to \bMM_{i,1}$ and
$q : \bMM_{i,1} \times \bMM_{g-i,1} \to \bMM_{g-i,1}$
are the projections to the first factor and the second factor
respectively. 
In the same way as for $\beta^*(D)$,
we can see that
$D_i$ (resp. $D_{g-i}$) is weakly positive at any points of $\MM_{i,1}$
(resp. $\MM_{g-i,1}$) by virtue of (1) of Corollary~\ref{cor:slope:cone:mg1}.
Thus, by using (2) of Proposition~\ref{prop:samp:pamp:pullback},
$\alpha_{i,g-i}^*(D)$ is weakly positive at any points
of $\MM_{i,1} \times \MM_{g-i,1}$. Therefore, we get (3) by
(2) of Proposition~\ref{prop:samp:pamp:push}.
\QED

\begin{Corollary}[$\ch(k) = 0$]
\label{cor:cone:mg:codim:two:another}
With notation as in Corollary~\rom{\ref{cor:cone:mg:codim:two}},
if $\rho_{irr}^*(D)$ is nef over $\rho_{irr}^{-1}(\bMM_g^{[1]})$
and $\rho_{i}^*(D)$ is nef over $\rho_{i}^{-1}(\bMM_g^{[1]})$ for all $i$,
then $D$ is nef over $\bMM_{g}^{[1]}$.
In particular, the Mori cone of $\bMM_g$ is the convex hull spanned
by curves lying on the boundary $\bMM_g \setminus \MM_g$,
which gives rise to a special case of \cite[Proposition~3.1]{GKM}.
\end{Corollary}

\Proof
Let $\beta' : \bMM_{g-1,2} \to \widetilde{\Delta}_{irr}$ and
$\alpha'_{i,g-i}:  \bMM_{i,1} \times \bMM_{g-i,1} \to 
\widetilde{\Delta}_{i}$ be the same as in Corollary~\ref{cor:cone:mg:codim:two}.
By our assumption,
$\beta^*(D) = {\beta'}^*(\rho_{irr}^*(D))$ is
nef over $\MM_{g-1,2}$ and
$\alpha_{i,g-i}^*(D) = {\alpha'}^*_{i,g-i}(\rho_{i}^*(D))$
is nef over $\MM_{i,1} \times \MM_{g-i,1}$ for every $i$.
Therefore, by Claim~\ref{thm:cone:mg:codim:two:claim:1}
in Theorem~\ref{thm:cone:mg:codim:two}, we can see that
$D$ is nef over $\bMM_{g}^{[1]}$.

Let $\Nef_{\Delta}(\bMM_g)$ be the dual cone of the convex hull
spanned by curves on the boundary $\Delta = \bMM_g \setminus \MM_g$.
In order to see the last assertion of this corollary, it
is sufficient to check $\Nef_{\Delta}(\bMM_g) = \Nef(\bMM_g)$,
which is a consequence of the first assertion.
\QED

\begin{Example}
For example,
the area of $(b_0, b_1)$ (resp. $(b_0, b_1, b_2)$) with
$\lambda - b_0\delta_0 - b_1\delta_1$ 
(resp. $\lambda - b_0\delta_0 - b_1\delta_1 -b_2 \delta_2$)
nef over $\bMM_3^{[1]}$ (resp.  $\bMM_4^{[1]}$)
is the inside of the following triangle (resp. polyhedron):
\ifneedpagebreakforpic
\[ \text{(see the next page)} \]
\vfill\eject
\else\fi
\[
\text{
\begin{pspicture}(-1,-2)(3,8.5)
\psline(0,-0.3)(2.08,-0.3)
\psline(0,-0.3)(2.78,7.2)
\psline(2.08,-0.3)(2.78,7.2)
\rput(0,-0.5){${\scriptscriptstyle (0,0)}$}
\rput(2.08,-0.5){${\scriptscriptstyle (\frac{1}{12},0)}$}
\rput(2.78,7.4){${\scriptscriptstyle (\frac{1}{9},\frac{1}{3})}$}
\rput(1.5,-1.5){{$\begin{array}{l}
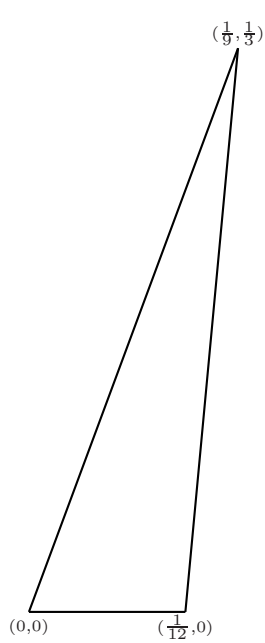
\text{{\tiny The area of $(b_0,b_1)$ with}} \\
\text{{\tiny $\lambda - b_0\delta_0 - b_1\delta_1$ nef over $\bMM_3^{[1]}$}}
\end{array}$}}
\end{pspicture}
}
\qquad\qquad\qquad
\text{
\begin{pspicture}(-1,-2)(4.25,8.5)
\psline(0,0)(0,3.4)
\psline(0,0)(1,-0.4)
\psline(0,0)(0.95,0.549) \psline(1.05,0.607)(2.25,1.3)
\psline(0,3.4)(1,3)
\psline(0,3.4)(4,7.4)
\psline(1,3)(4,7.4)
\psline(1,-0.4)(2.7,1.1)
\psline(2.25,1.3)(2.7,1.1)
\psline(2.7,1.1)(4,7.4)
\psline(2.25,1.3)(4,7.4)
\psline(1,-0.4)(1,3)
\rput(-0.5,0){${\scriptscriptstyle (0,0,0)}$}
\rput(1,-0.6){${\scriptscriptstyle (\frac{1}{12},0,0)}$}
\rput(3.3,1.1){${\scriptscriptstyle (\frac{1}{10},\frac{1}{5},0)}$}
\rput(1.7,1.45){${\scriptscriptstyle (\frac{1}{15},\frac{1}{5},0)}$}
\rput(4,7.65){${\scriptscriptstyle (\frac{1}{9},\frac{1}{3},\frac{4}{9})}$}
\rput(-0.5,3.4){${\scriptscriptstyle (0,0,\frac{2}{7})}$}
\rput(1.65,3){${\scriptscriptstyle (\frac{1}{12},0,\frac{2}{7})}$}
\rput(1.6,-1.5){{$\begin{array}{l}
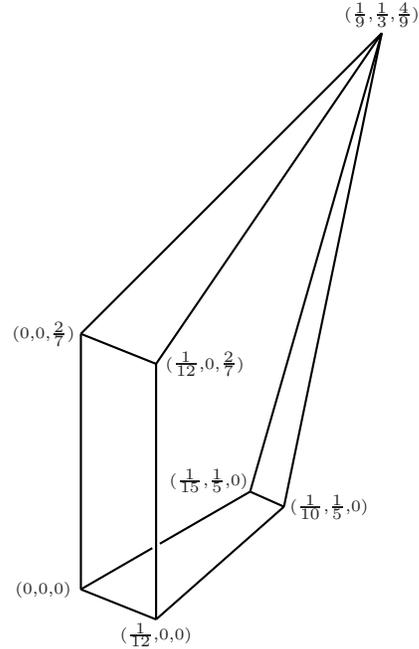
\text{{\tiny The area of $(b_0,b_1,b_2)$ with}} \\
\text{{\tiny $\lambda - b_0\delta_0 - b_1\delta_1 - b_2\delta_2$
nef over $\bMM_4^{[1]}$}}
\end{array}$}}
\end{pspicture}
}
\]
\end{Example}


\bigskip
\begin{thebibliography}{99}

\bibitem{CH}
M. Cornalba and J. Harris,
Divisor classes associated to families of stable varieties, with
application to the moduli space of curves, 
Ann. Scient. Ec. Norm. Sup., 21 (1988), 455--475.

\bibitem{DM}
P. Deligne and D. Mumford,
The irreducibility of the space of curves of given genus,
Publ. Math. IHES, 36 (1969), 75--110.

\bibitem{DeJ}
A. J. De Jong,
Smoothness, semi-stability and alteration,
Publ. Math. IHES, 83 (1996), 51--93.

\bibitem{FICM}
C. Faber,
Intersection-theoretical computation on $\bMM_g$,
Banach Center Publ. 36, Polish Acad. Sci., Warsaw, 71--81.
(alg-geom/9504005).

\bibitem{GKM}
A. Gibney, S. Keel and I. Morrison,
Towards the ample cone of $\bMM_{g,n}$, (math.AG/0006208).

\bibitem{Hain}
R. Hain,
Moriwaki's inequality and generalizations, talk at Chicago,
http://www.math.duke.edu/\~{}hain/talks/

\bibitem{KMCont}
S. Keel and J. McKernan,
Contractible extremal rays on $\bMM_{0,n}$,
preprint, (alg-geom/9607009).

\bibitem{KMBirat}
J. Koll\'{a}r and S. Mori,
Birational geometry of algebraic varieties,
Cambridge tracts in Mathematics, 134.

\bibitem{K23}
F. Knudsen,
The projectivity of the moduli space of stable curves, II and III,
Math. Scand. 52 (1983), 161--199 and 200--212.

\bibitem{Mo5}
A. Moriwaki,
A sharp slope inequality for general stable fibrations of curves,
J. reine angew. Math. 480 (1996), 177--195.

\bibitem{MoRB}
A. Moriwaki,
Relative Bogomolov's inequality and the cone of positive divisors
on the moduli space of stable curves,
J. of AMS, 11 (1998), 569--600.

\bibitem{MoCD}
A. Moriwaki,
The continuity of Deligne's pairing,
International Mathematics Research Notices, (1999), No.19, 1057--1066.

\bibitem{Se}
C. S. Seshadri,
Quotient spaces modulo reductive algebraic groups,
Ann. of Math., 95 (1972), 511--556.


\end{thebibliography}
\end{document}